\documentclass[12pt]{article}

\usepackage{amssymb,latexsym,amsmath, amsfonts}     
\usepackage{graphicx}
\usepackage{caption,subcaption}
\usepackage{bib entry}
\nobibliography*
\usepackage{algcompatible}
\usepackage{algorithm}
\usepackage{cleveref}
\usepackage{url}
\begin{document}

 \title{ \bf A bi-objective optimization framework for three-dimensional road alignment design}

\author{D. Hirpa\thanks{Combinatorics and Optimization, University of Waterloo, 200 University Avenue West
Waterloo, ON, N2L 3G1, Canada {\tt DessalegnH@aims.ac.za}}, W. Hare\thanks{Corresponding author: Mathematics, University of British Columbia, Okanagan Campus (UBCO), 3333 University Way, Kelowna, BC, V1V 1V7 Canada, {\tt Warren.Hare@ubc.ca}}, Y. Lucet\thanks{Computer Science, UBCO, {\tt Yves.Lucet@ubc.ca}},
Y. Pushak\thanks{Computer Science, UBCO, {\tt Yasha\underline{~}Pushak@hotmail.ca}}, S. Tesfamariam\thanks{Engineering, UBCO, {\tt Solomon.Tesfamariam@ubc.ca}}}
\maketitle

\begin{abstract}
Optimization of three-dimensional road alignments is   a nonlinear  non-convex optimization problem.  The development of models that fully optimize a three-dimensional road alignment problem is challenging due to numerous factors involved and complexities in the geometric specification of the alignment.  In this study, we developed a novel bi-objective optimization approach to solve a three dimensional road alignment problem where the horizontal and vertical alignments are optimized simultaneously. Two conflicting cost objective functions, \emph{earthwork} cost and the \emph{utility} cost, are cast in a bi-objective optimization problem. We numerically compare several multi-objective optimization solvers, and find that it is possible to determine the Pareto front in a reasonable time.
\end{abstract}


\section{Introduction}

Road design consists of finding a three-dimensional route alignment on the ground surface.
It aims to connect two terminals (the start and the end of a road) at minimum possible cost and maximal utility,  subject to  design, environmental, and social constraints \cite{jha2007multi}.  Since the number of alternative routes connecting two end points is unlimited, a traditional route location analysis, which has relied heavily on human judgement and intuition, may overlook many good alternatives  \cite{Chew1989315,li2013approach}. In order to consider all such route alternatives, and to reduce the workload for engineers,  several automated procedures that determine good road alignments and  approximate construction costs have been developed \cite{Jong2003107}. The automation of  the road design problem reduces the tedious and error-prone manual tasks, most notably drafting \cite{kim2004intersection}. In addition, this procedure allows the use of optimization techniques in search of a good alignment \cite{oecd1973}.  Optimization techniques save design time and  provide the decision maker with powerful tools that search for an alignment with minimum cost from a large number of alternative alignments. In fact, optimization of road alignment can yield considerable savings in construction costs when compared with unoptimized design procedures \cite{oecd1973}.

The road design problem  can be broken down into three interconnected stages: the  horizontal alignment, the vertical alignment, and the earthwork \cite{Hare2011}. The horizontal alignment is a bird's eye view of a road trajectory. A typical horizontal alignment is composed of a sequence of tangents, circular curves, and transition curves. Transition curves have the property that the radius of curvature changes progressively along them. The main considerations in horizontal alignment design are that it should avoid lands which are restricted or  expensive to purchase, obstacles which present engineering difficulties, and ground which may involve large amount of earthwork.  The cost of road construction for the horizontal alignment problem depends on the cost of acquiring land and on the output of the vertical alignment stage \cite{Hare2011}. Optimization of the horizontal alignment seeks a low cost route while adhering to the design standards and reducing environmental impacts  ~\cite{ahmad2011approach}.  However, optimization of horizontal alignment should also seek a highly utile route.  These two goals may often be in conflict with each other.  In the literature, the following  models have been developed for optimizing horizontal alignments: calculus of variation \cite{nicholson1973variational}, network optimization \cite{dan}, dynamic programming \cite{oecd1973}, and genetic algorithms \cite{preliminary}. Detailed discussion on the advantages and disadvantage of these methods can be found in \cite{schonfeld2006intelligent}.

The vertical alignment  is the view of the centreline of the road when seen along the longitudinal cross-section of the road.  A typical  vertical alignment is composed of  straight sections known as vertical tangents and  parabolic  curves, namely crest and sag curves. The design of   vertical alignment is a crucial step in the road design problem since it has implications on road construction costs, traffic operations, vehicle fuel consumption, and safety \cite{vertAligAnalysis}.  In  vertical alignment optimization, one fits a road profile to the ground profile while respecting various grade constraints and other road specifications. The objective is to minimize the cost of construction and the negative  impacts on the environment  \cite{vertAligAnalysis}.  Vertical alignment optimization models include linear programming \cite{easa1988selection, mayer1981earthmoving},  numerical  search \cite{hayman1970optimization}, state parametrization \cite{Goh1988399},  dynamic programming \cite{Goh1988399, goktepe2009optimiz}, genetic algorithm  \cite{ahmad2011approach, vertAligAnalysis}, and mixed integer linear programming \cite{Hare2014161, Hare2011}. (See also \cite{Goh1988399, goktepe2009optimiz, schonfeld2006intelligent} for more references).

A three-dimensional road alignment is the superimposition of two-dimensional horizontal and vertical alignments. In essence, road alignment design is a three-dimensional problem represented in $X$, $Y,$ and $Z$ coordinates.  The development of models that fully optimize a three-dimensional road design problem  is not yet successful, because there are many factors involved and complexities in the geometric specifications \cite{3dTabuSearch, goktepe2009optimiz, schonfeld2006intelligent}. Three-dimensional road alignment optimization is  presented as a constrained, nonlinear, and non-differentiable optimization problem \cite{3dTabuSearch}. There are two basic approaches found in the literature: models that simultaneously optimize the horizontal and vertical alignments \cite{3dTabuSearch,  Chew1989315, Jong2003107},  and models that employ two or more stages of optimization \cite{nicholson1973variational, parker3D}. Existing approaches for three-dimensional alignment optimization include genetic algorithms \cite{Jong2003107, Kang2012257, 3dAlign}, swarm artificial intelligence \cite{bosurgi2013pso}, dynamic programming \cite{li2013approach}, neighbourhood search \cite{3Dmodel}, and distance transform \cite{desmith}. Dynamic programming and numerical search methods suffer from the high computational effort and large memory requirements \cite{schonfeld2006intelligent, maji2012comparison, Jong2003107}.  So far, all of these method focus only on minimizing the overall design cost, and omit the conflicting cost of road utility.

This paper presents a novel bi-objective framework for optimizing a three-dimensional road alignment problem. In our model, a bi-objective cost function is formulated that consists of the `earthwork cost' and the `utility cost'. The cost associated with earthwork is the cost of ground cut, ground fill, waste material, and borrow material, all of which depend on the volume of the material.  The utility cost is computed based on the length of the road. Model constraints restrict maximum grade and radius of curvature for horizontal turns.

This paper advances research in horizontal alignment in several manners.  First, this paper considers the three-dimensional alignment as a bi-objective problem. Research in this direction has been limited, but shown great promise in application  \cite{jha2007multi, ATR5670430405, Maji2011966}.  Indeed, aside from \cite{jha2007multi, ATR5670430405, Maji2011966} (which all present different aspects of the same model), research in this direction appears absent.  The principal difference in our model, as compared to \cite{jha2007multi, ATR5670430405, Maji2011966}, is expanded details on exactly how the road alignment is designed and costs are computed.   As a result, our model also provides a detailed new surrogate model that can produce both an earthwork cost and a utility costs for a potential alignment.  This is a very significant contribution to the field, as it makes future research in this area more accessible.  We show that the model is computationally tractable, but retains many of the key aspects required to model earthwork and utility costs.  Finally, a case study is provided that demonstrates the Pareto front of the bi-objective problem can be computed in reasonable time.  This case study compares three different multiobjective optimization algorithms namely, the multiobjective genetic algorithm (GAMOO) \cite{deb2001multi}, the direct multisearch for multiobjective optimization (DMS) \cite{dmsMOO}, and the weighted sum method (WS) \cite{normalization, weightedsum}.

The remainder of this paper is organized as follows.  In  Section \ref{modelvar}, the  models for both horizontal and vertical alignments are presented. The  cost functions are formulated in Section \ref{surrcost}. The bi-objective optimization model is  provided in  Section \ref{mop}.  The numerical experiments performed  and their results are reported  in Section \ref{numexp}.    Finally, Section \ref{colcln} presents some concluding remarks.


\section{ Model variables }\label{modelvar}

\subsection{ The geometric design for horizontal alignment}
\begin{figure}[ht]
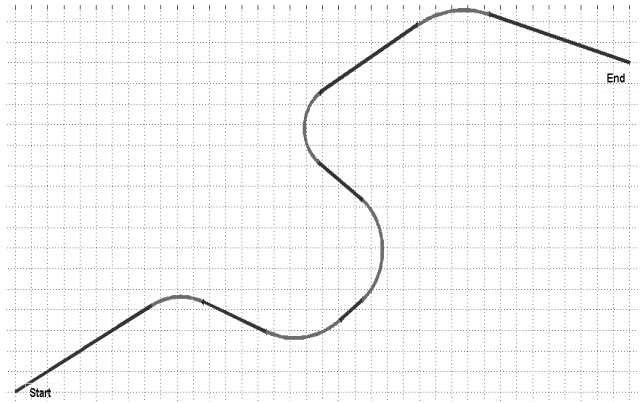
FIGURE 1 HERE\caption{An example of  horizontal alignment.}\label{fig11}\end{figure}

In our  model,  a horizontal alignment is composed of  tangential segments and circular curves. The horizontal alignment geometry is required to satisfy two criteria: (1) the alignment should satisfy the  orientation specified as $tangent-circle-tangent$, and (2) the first and the last sections of the alignment should be tangent segments.  The circular curves are placed between two adjacent tangents to mitigate the sudden changes in the direction of the alignment, as is standard in road design.  A typical geometry of horizontal alignment is depicted in Figure \ref{fig11}.

The geometric design of a horizontal alignment is described by a sequence of
intersection points $IP$ and  the radius of a circular curve $r$.   There are three decision variables associated with each intersection point, namely  $x,$  $y,$ and $r$. The   $x$ and  $y$ variables correspond to  the $x-$ and $y-$coordinates of the intersection point, and $r$ is the radius of the circular curve.

\begin{figure}[ht]
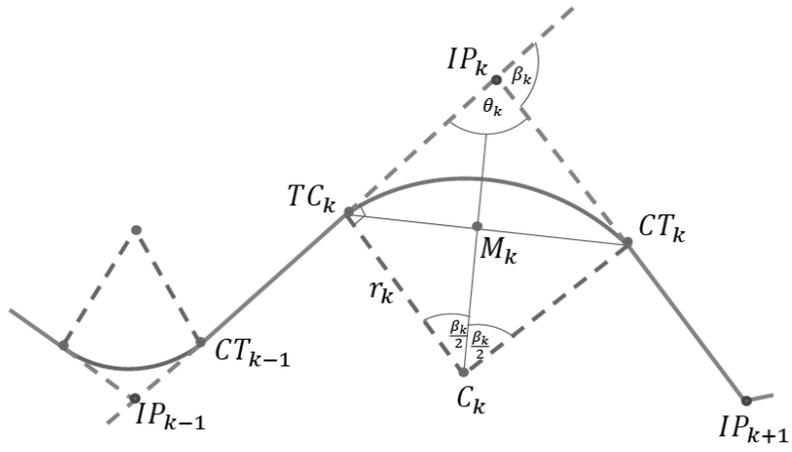
FIGURE 2 HERE\caption{Horizontal alignment geometry}\label{figg}\end{figure}

Let $START=(x_{s},y_{s},z_{s})$ and $END=(x_{e},y_{e},z_{e})$ be the start and end  of the design road. We denote the horizontal component of  $START$ and $END$ by  $S_{h}=(x_{s},y_{s})$ and $E_{h}=(x_{e},y_{e}).$  Given the set of $N$ intersection  points
            \[ \lbrace IP_{1},IP_{2},\cdots, IP_{N}\rbrace ,\]
and the corresponding set of radius of curvature
			\[ \lbrace r_{1},r_{2},\cdots, r_{N}\rbrace ,\]
one can insert $ N + 1$ tangential road segments and $N$ circular curves to get the exact shape of a horizontal alignment (see Figure \ref{figg}).
Next, we calculate the center $C_{k}$ of the curve and transition points $TC_k$  and $CT_k$, as denoted in Figure \ref{figg}.
 Denote  the vector from $IP_{k}$ to $IP_{k-1}$ and $IP_{k}$ to $IP_{k+1}$, respectively, by ${\bf IP}_{k(k-1)}$ and ${\bf IP}_{k(k+1)}$.
  Define $C_{k} =(x_{c_{k}},y_{c_{k}})$, $TC_{k}=(x_{tc_{k}},y_{tc_{k}})$ and $CT_{k}=(x_{ct_{k}},y_{ct_{k}})$.  Let $\theta_{k}$ be the  angle between ${\bf IP}_{k(k-1)}$ and ${\bf IP}_{k(k+1)}$, then
 we calculate the coordinates of $C_{k}$, $TC_{k}$, and $CT_{k}$ as follows.

\subsubsection{Calculation of the transition points $TC_{k}$ and $CT_{k}$}\label{trans}
Given the intersection points $IP_{k-1},IP_{k}$ and $IP_{k+1}$, and radius $r_{k}$, see Figure \ref{figg},  then

\begin{equation}\label{angle}
\theta_{k}=\arccos\bigg(\frac{{\bf IP}_{k(k-1)}\cdot{\bf IP}_{k(k+1)} }{\Vert{\bf IP}_{k(k-1)} \Vert \Vert {\bf IP}_{k(k+1)}\Vert}\bigg),
\end{equation}
where ${\bf IP}_{k(k-1)}=IP_{k-1}-IP_{k}~\text{and}~
{\bf IP}_{k(k+1)}=IP_{k+1}-IP_{k}$.
 Let $L_{k}$ be the distances from $IP_{k}$ to $TC_{k}$.
Then
\begin{align*}
 L_{k}&=r_{k}\tan\frac{\beta_{k}}{2},
\end{align*}
where $\beta_{k}$ is the central angle of the circular arc, see Figure \ref{figg}.

Using the identity
\[\tan\frac{\beta_{k}}{2}= \frac{1-\cos\beta_{k}}{\sin\beta_{k}}~\text{ and the fact that}~\beta_{k}=\pi-\theta_{k},\]
we have
\[L_{k}=r_{k}\dfrac{1+\cos\theta_{k}}{\sin\theta_{k}}.\]
Note that the distance from  $IP_{k}$ to $CT_{k}$ is  equal to $L_{k}.$

Define vectors ${\bf V_{k}}$ and ${\bf U_{k}}$ as
\[{\bf V}_{k}=TC_{k}-IP_{k}~ \text{and}~{\bf U}_{k}=CT_{k}-IP_{k}.\]
 Then
\[{\bf V}_{k}= TC_{k}-IP_{k}=L_{k}\dfrac{{\bf IP}_{k(k-1)}}{\Vert {\bf IP}_{k(k-1)}\Vert},~ \text{and}~{\bf U}_{k}= CT_{k}- IP_{k}=L_{k}\dfrac{{\bf IP}_{k(k+1)}}{\Vert {\bf IP}_{k(k+1)}\Vert}. \]
Thus, we obtain
\begin{equation}\label{tangent}
TC _{k}= IP_{k}+L_{k}\dfrac{{\bf IP}_{k(k-1)}}{\Vert {\bf IP}_{k(k-1)}\Vert}, ~\text{and}~
 CT _{k}= IP_{k}+L_{k}\dfrac{{\bf IP}_{k(k+1)}}{\Vert {\bf IP}_{k(k+1)}\Vert}.
\end{equation}

\subsubsection{Calculation of the center $C_{k}$}

Let $M_{k}$ be the midpoint of the line segment that connects $TC_{k}$ and $CT_{k}$. Clearly $M_{k}$ lies on the line that passes through $IP_{k}$ and $C_{k}$. Hence
\[M_{k}=\dfrac{1}{2}\big(TC_{k}+ CT_{k}\big)= IP_{k}+\dfrac{L_{k}}{2}\bigg( \dfrac{{\bf IP}_{k(k-1)}}{\Vert {\bf IP}_{k(k-1)} \Vert}+\dfrac{{\bf IP}_{k(k+1)}}{\Vert {\bf IP}_{k(k+1)} \Vert}\bigg).\]
If  $d_{k}$ is  the distance from $IP_{k}$ to the circle that is  calculated as $$ d_{k}=r_{k}\sec\frac{\beta_{k}}{2}-r_{k}=r_{k}\csc\frac{\theta_{k}}{2}-r_{k},~\text{since}~ \sec\frac{\beta_{k}}{2}=\dfrac{r_{k}+d_{k}}{r_{k}},$$
then
\begin{equation}\label{centre}
 C_{k}- IP_{k}=(r_{k}+d_{k})\dfrac{ M_{k}- IP_{k}}{\Vert  M_{k}- IP_{k}\Vert} ~ \text{so}~ C_{k}= IP_{k}+r_{k}\csc\frac{\theta_{k}}{2}\bigg(\dfrac{ M_{k}-IP_{k}}{\Vert  M_{k}- IP_{k}\Vert}\bigg).
\end{equation}

\subsection{Vertical alignment  geometric design}
\begin{figure}[ht]
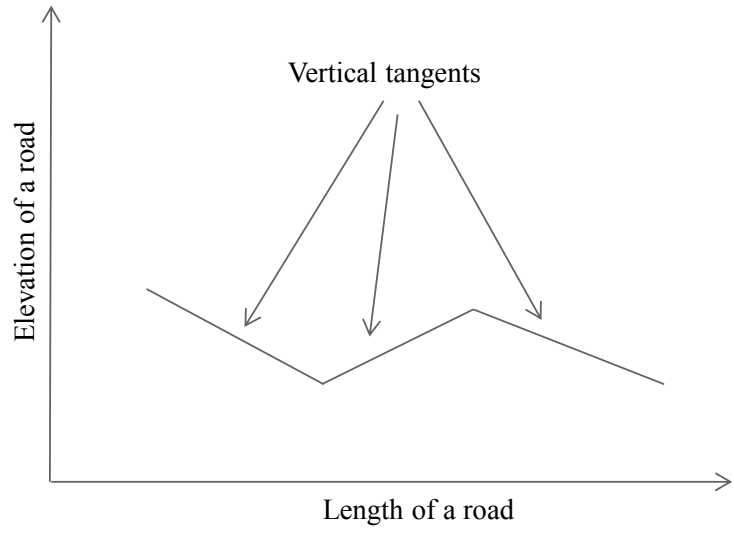
FIGURE 3 HERE\caption{An example of vertical alignment model.}\label{visina8}\end{figure}

\begin{figure}[ht]
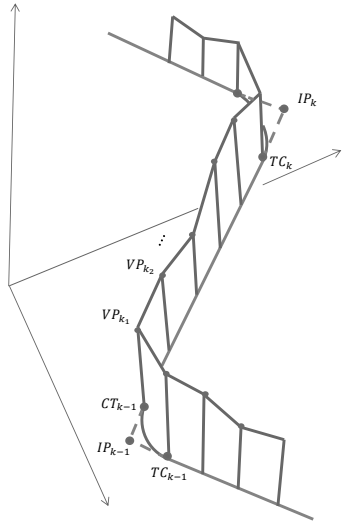
FIGURE 4 HERE\caption{Road alignment in 3D.}\label{fig321}\end{figure}

The vertical alignment geometric design specifies the elevation of a point along
 a roadway. In our model it is composed of straight lines known as vertical tangents.   We begin the model by stretching  a surface orthogonal to the $xy-$plane along the horizontal alignment. The flattened surface is called   the $hz-$plane, where $h$ is the distance measured along the horizontal alignment. The projection onto the $hz-$plane of the three-dimensional road alignment is the vertical alignment \cite[page~36]{schonfeld2006intelligent}. A typical vertical alignment model is shown in Figure \ref{visina8}.
 Given the  following sets of intersection points and the horizontal tangent points, calculated in Section \ref{trans},
\begin{align}
{\bf IP}&= \lbrace  IP_{1},IP_{2},\cdots, IP_{N}\rbrace,\\
{\bf TC}&=\lbrace  TC_{1},TC_{2},\cdots, TC_{N}\rbrace, and\\
{\bf CT}&=\lbrace  CT_{1},CT_{2},\cdots, CT_{N}\rbrace,
\end{align}
we  partition the horizontal tangent between  $CT_{k-1} ~\text{and}~ TC_{k}$ into $m_{k}$ equally spaced points, and  corresponding to these points we define a set of vertical points
 \begin{align}\label{vp}
\textbf{VP}_{k}=& \lbrace VP_{k,1},VP_{k,2}, \cdots, VP_{k,m_{k}}\rbrace,
\end{align}
where $VP_{k,j}=(x_{k,j},y_{k,j},z_{k,j})$, see Figure \ref{fig321}.

The mathematical details of the design road elevation $z_{k_{j}}$, $x_{k_{j}}$ and $y_{k_{j}}$ is given below.
We require the elevation to be linear over circular road sections, which reduces the computation time. This assumption does not have significant effect on the precision of the solution since the length of the curved section is typically much smaller than the tangent section.

Define $p_{j-1}=(x_{k,j-1},y_{k,j-1})$ and $p_{j}=(x,_{k,j},y_{k,j})$.  Then,  for  $j=1,2,\cdots,m_{k}$
 \begin{subequations}
\begin{align}
 x_{k,j}&=x_{ct_{k-1}}+\dfrac{j}{m_{k}}(x_{tc_{k}}-x_{ct_{k-1}}),~~~~~ ~~
 y_{k,j}=y_{ct_{k-1}}+\dfrac{j}{m_{k}}(y_{tc_{k}}-y_{ct_{k-1}}), \label{vpts2}\\
 x_{k,j-1}&=x_{ct_{k-1}}+\dfrac{j-1}{m_{k}}(x_{tc_{k}}-x_{ct_{k-1}}),~~
 y_{k,j-1}=y_{ct_{k-1}}+\dfrac{j-1}{m_{k}}(y_{tc_{k}}-y_{ct_{k-1}})\label{vpts4}.
\end{align}
\end{subequations}

Note that $z_{k,j-1}$ and $z_{k,j}$ are  decision variables for the vertical alignment corresponding to   $p_{j-1}$ and $p_{j}$. Now we parametrize the  vertical tangent between $VP_{k,j-1}$ and $VP_{k,j}$ using the parameter $s$ as follows. For $j=1,2,\cdots,m_{k}$,  the coordinates of any point $p=(x,y,z)$ on the line segment between $VP_{k,j-1}$ and $VP_{k,j}$ are calculated by
\begin{subequations}
\begin{align}
x(s)&=x_{k,j-1}+(x_{k,j}-x_{k,j-1})(m_{k}s+1-j)=x_{ct_{k-1}}+s(x_{tc_{k}}-x_{ct_{k-1}}), \label{param1} \\
y(s)&=y_{k,j-1}+(y_{k,j}-y_{k,j-1})(m_{k}s+1-j)=y_{ct_{k-1}}+s(y_{tc_{k}}-y_{ct_{k-1}}), \label{param2} \\
z(s)&=z_{k,j-1}+(z_{k,j}-z_{k,j-1})(m_{k}s+1-j). \label{param3}
\end{align}
\end{subequations}
where $\dfrac{j-1}{m_{k}}\leq s\leq \dfrac{j}{m_{k}}.$

The horizontal distance between $VP_{k,j-1}$ and $VP_{k,j}$, denoted  $d_{k,j}$, is calculated as
\[d_{k,j}=\sqrt{(x_{k,j}-x_{k,j-1})^{2}+(y_{k,j}-y_{k,j-1})^{2}}.\]
  The number of vertical alignment  variables  on each horizontal tangent section, except for the first and last horizontal tangents, is $m_{k}+1$. The first and last horizontal tangent section each have  $m_{0}$ and $m_{N}$ variables, respectively,  because the start and the end of a  road are not variables.
Suppose the number of  intersection points is $N$. Then the number of vertical alignment variables  is
\[M=m_{0}+m_{N}+\sum_{k=1}^{N-1}(m_{k}+1).\]


\subsection{Model Constraints}\label{const}
The geometric design of horizontal alignments and vertical alignments are constrained by code requirements and other limitations. In this section,  the set of constraints for both alignments are discussed.

\subsubsection*{Box constraint on the IP location}
 Let $IP_{k}=(x_{k},y_{k}),$ and let $x_{u,k},x_{l,k},y_{u,k},~\text{and}~y_{l,k}$ be real numbers. Then the box constraint corresponding to $IP_{k}$ is given by
\[x_{l,k}\leq x_{k}\leq x_{u,k},~\text{and}~ y_{l,k}\leq y_{k}\leq y_{u,k}, ~~k=1,2,\cdots,N.\]
These constraints help reduce the search space for potential designs, and thereby allow for reasonable sized data sets (see Section \ref{surrcost}.)
\subsubsection*{The circular curves should not overlap}
In our model, two adjacent circular curves portions of a road are allowed to meet, but never overlap. Overlapping would result in a discontinuous road. This requirement can be written mathematically as follows
\begin{equation}\label{Hconst1}
0\leq \Vert TC_{k}-IP_{k-1}\Vert -\Vert CT_{k-1}-IP_{k-1}\Vert, ~~k=1,2,\cdots,N.
\end{equation}
\subsubsection*{The minimum turning radius}
Given the minimum radius $r_{\text{min}}$, the  radius of curvature at each intersection point has to satisfy the minimal radius requirement
\begin{equation}
r_{\text{min}}\leq r_{k}, ~~k=1,2,\cdots,N.
\end{equation}
\subsubsection*{Maximum allowable slope or gradient of the road}
 Given the maximum allowed vertical slope, $G_{\max}$,
\[\vert z_{k,j}-z_{k,j-1}\vert\leq d_{k,j}G_{\max},~\text{for all }~k~\text{and}~j.\]
\subsubsection*{The design road elevation bounds}
While not strictly necessary, it is reasonable to bound the road elevation.
  Given the maximal vertical offset $\bar{z}$ at the current ground elevation $z_{g}$,
  \[z_{g}-\bar{z}\leq z_{k,j}\leq z_{g}+\bar{z},~\text{for all }~k~\text{and}~j.\]
\subsubsection*{Circular road sections are linear}
  In our model, we require the elevation to be linear over circular road sections. This simply means that the circular road sections are not split into multiple vertical sections. 


\section{Surrogate cost formulation}\label{surrcost}
\begin{figure}[ht]
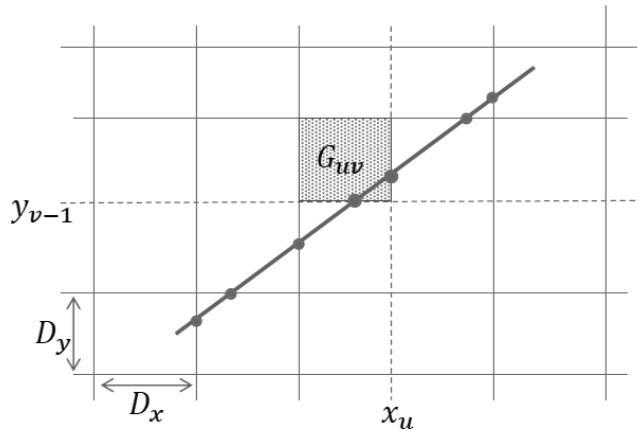
FIGURE 5 HERE\caption{An example, the projection of a tangent road segment onto the horizontal plane.}\label{fig2d}
\end{figure}

The  surrogate cost model  developed is based on the approximation of the ground surface by a series of continuous  linear functions or planar surfaces. The horizontal plane in the region of interest is divided into grid cells of equal size, which are small enough that the terrain above the grid is approximated by a single linear function. A matrix format is employed to store important ground information for the  region of interest.  We denote the $(u,v)$  grid cell by $G_{uv}$, see Figure \ref{fig2d}. The terrain in each grid is then approximated by  a linear function
\[ z=A_{uv}x+B_{uv}y+C_{uv},\]
where $x,y\in G_{uv},$ and $A, B$ and $C$ are matrices containing the ground elevation data.   Hereafter, we refer to the approximating plane elevation as the ground profile.
We aim to calculate the cost of a road section (tangential section or circular section), which is the sum of costs corresponding to each grid cell along the section.

In each section, two costs are developed: the earthwork costs and the utility cost. The cost associated with earthwork is the cost of ground cut, ground fill, and cost of borrow or waste material, all of which depend on the volume of the material. Borrow and waste material is the difference between ground cut and ground fill. The calculation of the volume of ground cut and ground fill depends on the elevation difference between the approximating plane and the design road. The utility cost is computed based on the length of the road.

For the collection $S$ of tangential road sections and circular road sections, we define $V_{c_{\xi}}$ as the total volume of ground cut on road section $\xi\in S$, and $V_{f_{\xi}}$ as the total volume of ground fill on road section $\xi\in S$. The overall volume of cut, denoted $V_{c}$ and volume of fill, denoted $V_{f}$, is given as
\[
V_{c}=\sum_{\xi\in S}V_{c_{\xi}},~\text{ and}~
V_{f}=\sum_{\xi\in S}V_{f_{\xi}}.
\]
We then apply the costing parameters to calculate the total cost from earthwork $Cost_{\tt e}$ and the total utility cost $Cost_{\tt u}$,
\begin{equation}\label{surrogateCost}
Cost_{\tt e}= C_{c}V_{c}+C_{f}V_{f}+C_{w}\Vert V_{f}-V_{c}\Vert,~ Cost_{\tt u}=C_{u}L ,
\end{equation}
where $C_c, C_f,$ are parameters for cut and fill costs, $C_w,$ is a single parameter representing the cost of either waste or borrow, $C_u$ is a utility cost parameter, and $L$ is the length of the road.  Note that waste and borrow costs are combined into a single cost, which represents the amount of unbalanced material in the road design.  Also note that $C_c$, $C_f$, and $C_w$ can be viewed as cost parameters in present day dollar (the `construction dollars').  Conversely, $C_u$ should be viewed as a parameter representing `road effectiveness', so does not necessarily have a clear translation to present day dollars.

In the following sections, the calculation of cost models for tangent road section and circular road section are provided.

\subsection{Surrogate cost model for tangent road segment}
\label{sec:Surrogate cost model for tangent road segment}
Suppose $\xi_{k} \in S$ corresponds to a tangential road section that connects $VP_{k,j-1}$ to $VP_{k,j}$. The parametric equation of $\xi_{k}$ is given by
\[r_{t}(s) =\big(x(s), y(s),z(s)\big),\]
where  $x(s),y(s)$, and $z(s)$ are given in \cref{param1,param2,param3}.

The parametric equations of the design road profile and the  ground profile in grid $G_{uv}$ is calculated  as
\begin{equation}\label{RdElev}
z_{r}(s)=z_{k,j-1}+(z_{k,j}-z_{k,j-1})(m_{k}s+1-j),
\end{equation}
and
\begin{align}\label{GrdElev}
z_{g}(s)&=A_{uv}x(s)+B_{uv}y(s)+C_{uv} \nonumber\\
        &=\big( A_{uv}x_{ct_{k-1}}+ B_{uv}y_{ct_{k-1}}+C_{uv}\big)+
        s\big( A_{uv}(x_{tc_{k}}-x_{ct_{k-1}}) + B_{uv}(y_{tc_{k}}-y_{ct_{k-1}})\big),
\end{align}
where $x(s),y(s)\in G_{uv}.$
\subsubsection{ Parameter calculation}
We compute the  parameter $s$ corresponding to the $x$ boundary cross and the $y$ boundary cross of  grid cell $G_{uv}$. Since it may be required to  cut and fill the ground within a particular grid cell, we also  calculate the parameter corresponding to  the point of transition from cut to fill or fill to cut. These parameters are used to  calculate the design road elevation in \eqref{RdElev} and the ground elevation in \eqref{GrdElev}.
We define the  sets of parameters $T_{x_{k}}^{j},T_{y_{k}}^{j},T_{t_{k}}^{j}$  for the $x$-boundary cross, $y$-boundary cross, and for the   transition point.

Let $x_{u},$ and $ y_{v-1}$ be the $x$ and $y$ boundary points  at which the horizontal tangent   crosses the grid cell $G_{uv},$  see Figure \ref{fig2d}. Then, for $\frac{j-1}{m_{k}}\leq s \leq \frac{j}{m_{k}}$ the  parameter $s$ corresponding to the boundary point  is given as
 \[T_{x_{k}}^{j}= \big \lbrace s~|~x(s)=x_{u} \big \rbrace ~\text{and}~ T_{y_{k}}^{j}=\big \lbrace s~|~ y(s)=y_{v-1}\big \rbrace.\]
Similarly, if there exists a  transition from cut to fill (or fill to cut) within $G_{uv}$,  the  transition set is calculated by
\[ T_{t_{k}}^{j}=\big\lbrace s~| ~z_{g}(s)=z_{r}(s)\big\rbrace.\]

 Define the union of the sets of parameters as
\[T_{k}^{j}=T_{x_{k}}^{j}\cup T_{y_{k}}^{j}\cup T_{t_{k}}^{j}.\]
Let $K$ be the cardinality of $T_{k}^{j}$. We sort $T_{k}^{j}$ in an increasing order to create $T_{k,s}^{j}$ as
\[T_{k,s}^{j}=\lbrace s_{1}<\cdots<s_{K-1}<s_{K}\rbrace. \]
\subsubsection{ Length of the tangent road section}

Given the vertical points $VP_{k,j}=(x_{k,j},y_{k,j},z_{k,j})$ and $VP_{k,j-1}=(x_{k,j-1},y_{k,j-1},z_{k,j-1})$, where $x_{k,j},x_{k,j-1},y_{k,j}$, and $y_{k,j-1}$ are calculated using   \cref{vpts2,vpts4}.
The length of $r_{t}(s)$, denoted $L_{t_{k}}^{j}$, is calculated by
\begin{align}\label{des}
L_{t_{k}}^{j}=\dfrac{1}{m_{k}^{2}}\sqrt{(x_{tc_{k}}-x_{ct_{k-1}})^{2}+(y_{tc_{k}}-y_{ct_{k-1}})^{2}+ m_{k}^{2}(z_{k,j}-z_{k,j-1})^{2}}.
\end{align}
The length of a tangential road segment is computed as
\begin{equation}
L_{t_{k}}=\sum_{j=1}^{m_{k}}L_{t_{k}}^{j}.
\end{equation}

\subsubsection{ Volume of ground cut}

  If the parameters $s_{i-1},s_{i} \in T_{k,s}^{j}$ bracket a cut region, then the approximate volume of  ground cut   for the tangential road segment over $[s_{i-1},s_{i}]$ is obtained by integrating the  cross-section area.   The elevation difference between the design road and the ground at $s\in [s_{i-1},s_{i}]$   is given by
\begin{align*}\label{elev}
h_{c}(s)= z_{g}(s)-z_{r}(s)
&=\big(  A_{uv}x_{ct_{k-1}}+ B_{uv}y_{ct_{k-1}}- z_{k,j}+j(z_{k,j}-z_{k,j-1})+C_{uv}\big) +\\
&~~~~s\big(A_{uv}(x_{tc_{k}}-x_{ct_{k-1}}) + B_{uv}(y_{tc_{k}}-y_{ct_{k-1}})-m_{k}(z_{k,j}-z_{k,j-1}) \big).
\end{align*}

\begin{figure}[ht]
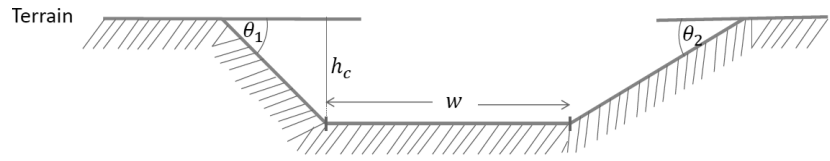
FIGURE 6 HERE\caption{An example of cut cross-section. }\label{fig24}\end{figure}

We assume that the ground cut has a trapezoidal cross-section, see Figure \ref{fig24}.

The cross-sectional area of the trapezoid  is calculated by
\begin{align}
a_{c}&=\frac{1}{2}h_{c}^{2}\cot\theta_{1}+Wh_{c}+\frac{1}{2}h_{c}^{2}\cot\theta_{2} =h_{c}\big(W+\frac{1}{2}h_{c} \kappa \big)
\end{align}
where  $W$ is the width of the road, $\theta_{1}$ and $\theta_{2}$ are side slope angles, and  $ \kappa = \cot\theta_{1}+\cot\theta_{2}$.

   Typically, for the model, we fix the values of $\theta_{1}$ and $\theta_{2}$  for all road segments; however, these values can change from segment to segment if desired.  In fact, these values can be based on the particular grid cell ($G_uv$) through which the road travels.  As such, local soil conditions can be incorporated into the model if desired.
The approximate volume of ground cut over $[s_{i-1},s_{i}]$ is calculated as
 \begin{equation}\label{5}
 V_{c_{k}}^{j}=\int_{s_{i-1}}^{s_{i}}a_{c}(s)||r_t'(s)||ds,
 \end{equation}
 where $j=1,2,3\cdots,m_{k}.$

Assuming $s_{i-1}$ and $s_i$ bracket a cut, the cross-section area $a_c(s)$ is given by
\begin{align}\label{area11}
a_{c}(s)&= Wh_{c}(s)+    \frac{1}{2}\kappa h_{c}^{2}(s)
 =W(\Omega +s\Theta) +\frac{1}{2}\kappa (\Omega +s\Theta)^{2}\nonumber \\
 &=\big( W\Omega+\frac{1}{2}\kappa \Omega^{2} \big) + \big(  W\Theta+\kappa\Omega\Theta \big)s + \frac{1}{2}\kappa\Theta^{2}s^{2},
\end{align}
where
\begin{align}
\Omega &=\big( A_{uv}x_{ct_{k-1}}+ B_{uv}y_{ct_{k-1}}- z_{k,j}+j(z_{k,j}-z_{k,j-1})+C_{uv} \big) \nonumber \\
 \Theta &= \big(A_{uv}(x_{tc_{k}}-x_{ct_{k-1}}) + B_{uv}(y_{tc_{k}}-y_{ct_{k-1}})-m_{k}(z_{k,j}-z_{k,j-1})\big).
\end{align}
If $s_{i-1}$, $s_{i}$ do not bracket a cut, then $a_c(s)$ is defined as $0$. We have $s_{i-1},s_{i} \in T_{k,s}^{j}$, so we can calculate
\begin{align}\label{Tangent Jacobian}
||r_t'(s)|| = ||r_t(s_{i}) - r_t(s_{i-1})||/(s_i - s_{i-1}),
\end{align}
since $r_t(s)$ is simply a line connecting $r_t(s_{i-1})$ and $r_t(s_i)$
over $[s_{i-1},s_{i}]$.
 Hence, if $s_{i-1}, s_{i}$ brackets a cut segment, integrating \eqref{5} using \eqref{area11} we get
\begin{align*}
V_{cut}=(\big( W\Omega+\dfrac{1}{2}\kappa \Omega^{2} \big)(s_{i}-s_{i-1})+
\dfrac{1}{2}\big(  W\Theta+\kappa\Omega\Theta \big)(s_{i}^{2}-s_{i-1}^{2})+\\ \frac{1}{6}\kappa\Theta^{2}(s_{i}^{3}-s_{i-1}^{3}))(||r_t(s_{i}) - r_t(s_{i-1})||/(s_i - s_{i-1})).
\end{align*}
Therefore, the volume of ground cut is computed as
\begin{align}\label{volumeCut}
V_{c_k}^j = \left\{ \begin{array}{rl}
 	V_{cut}~ & \mbox{if}~ s_{i-1}, s_i~ \mbox{bracket a cut segment} \\
	0 & \mbox{otherwise} .
\end{array}\right.
\end{align}


\subsubsection{ Volume of ground fill}
\begin{figure}[ht]
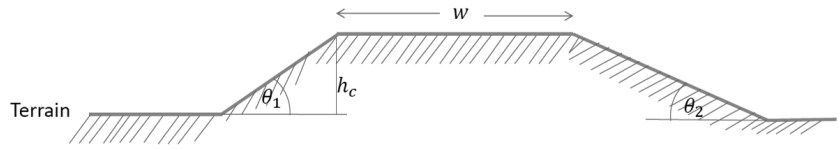
FIGURE 7 HERE \caption{An example, fill cross-section. }\label{fig25}\end{figure}

 Suppose the parameters $s_{i-1}$ and $s_{i}$ bracket a fill region. The elevation difference between the design road and the ground at  $s\in [s_{i-1},s_{i}]$ is given as
 \[h_{f}(s)=-h_{c}(s).\]
 We also assume that, the ground fill has trapezoidal cross-section, see Figure \ref{fig25}. If $s_{i-1}$ and $s_{i}$ bracket fill, then  the area of fill volume, $a_{f}$, is given as
 \begin{align}\label{area1}
a_{f}(s)&= Wh_{f}(s)+    \frac{1}{2}\kappa h_{f}^{2}(s)
=W(\Omega +s\Theta) +\frac{1}{2}\kappa (\Omega +s\Theta)^{2}\nonumber \\
 &=\big( W\Omega+\frac{1}{2}\kappa \Omega^{2} \big) + \big(  W\Theta+\kappa\Omega\Theta \big)s + \frac{1}{2}\kappa\Theta^{2}s^{2},
\end{align}
where
\begin{align}
\Omega &=\big( z_{k,j}-A_{uv}x_{ct_{k-1}}- B_{uv}y_{ct_{k-1}}-j(z_{k,j}-z_{k,j-1})-C_{uv} \big) \nonumber \\
 \Theta &= \big(m_{k}(z_{k,j}-z_{k,j-1})-A_{uv}(x_{tc_{k}}-x_{ct_{k-1}}) - B_{uv}(y_{tc_{k}}-y_{ct_{k-1}})\big).
\end{align}
If $s_{i-1}$ and $s_{i}$ do not bracket fill, then $a_{f}$ is calculated  as zero. Therefore, the volume of ground fill is calculated as
 \begin{align}
V_{f_k}^j = \left\{ \begin{array}{rl}
 	\textstyle \int_{s_{i-1}}^{s_{i}}a_{f}(s)||r_t'(s)||ds ~ & \mbox{if}~ s_{i-1}, s_i~ \mbox{bracket a fill segment} \\
	0 & \mbox{otherwise} .
\end{array}\right.
\end{align}


\subsection{ Surrogate cost model for circular road section}
\label{sec:Surrogate cost model for circular road section}
\begin{figure}[ht]
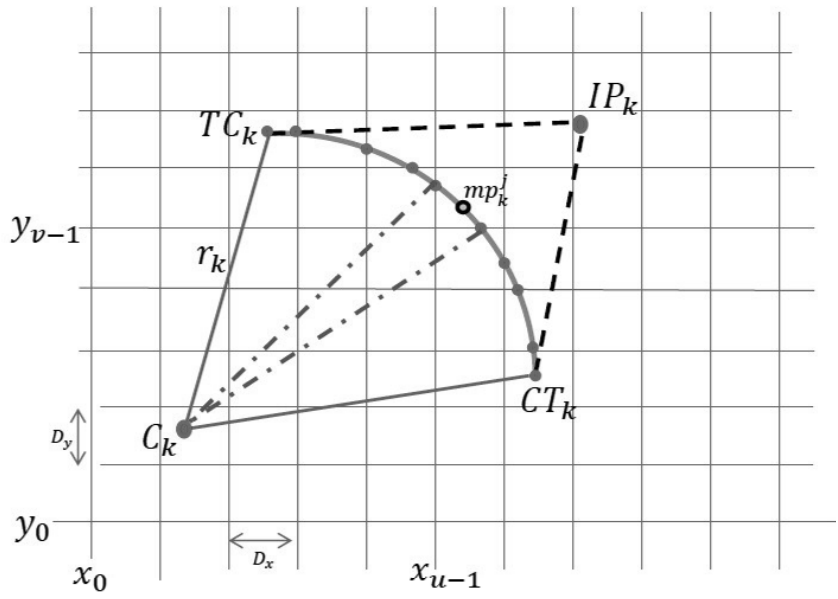
FIGURE 8 HERE\caption{Example of a horizontal curve section}\label{fig27}\end{figure}

Suppose  segment $\xi_{k} \in S$ corresponds to a circular road section, that is, a circle of radius $r_{k}$ and centre $(x_{c_{k}},y_{c_{k}})$ that connects $TC_{k}$ and $CT_{k}$, see Figure \ref{fig27}.  The parametric equation of $\xi_{k}$ is given as
\begin{equation}
r_{c}(s)=(x(s),y(s),z(s)),
\end{equation}
where,
\begin{align}
x(s)&=x_{c_{k}}+r_{k}\cos\theta_{c_{k}}(s), \nonumber \\
y(s)&=y_{c_{k}}+r_{k}\sin\theta_{c_{k}}(s),\nonumber \\
z(s)&=z_{k,m_{k}}+(z_{(k+1),1}-z_{k,m_{k}})s,~\text{and}\nonumber \\
\theta_{c_{k}}(s)&=\theta_{TC_{k}}+(\theta_{CT_{k}}-\theta_{TC_{k}})s.\nonumber
\end{align}

Similar to tangent segments, we begin by creating the parameter collection representing when the segment crosses an x or y boundary, or changes from a cut to fill;
\[T_{x}^{k}=\lbrace s~|~x(s)=x_{u-1}\rbrace,~~  T_{y}^{k}=\lbrace s~|~y(s)=  y_{v-1} \rbrace, ~\text{and}~T_{t}^{k}=\lbrace s~|~z_{g}(s)=z_{r}(s)\rbrace .\]
We then order them to create
 \begin{align}
T_{k,s}&=\lbrace 0=s_{1}< s_{2}<s_{3}<\cdots,s_{l-1}<s_{l}=1\rbrace. \label{paramSorted}
\end{align}

\subsubsection{ Length of  circular road section}

The length of the  circular road section, denoted by $L_{c_{k}}$ is given by
\begin{equation}
L_{c_{k}}=\sqrt{(r_{k}\theta_{c_{k}})^2 + (z_{(k+1),1}-z_{k,m_{k}})^2},
\end{equation}
where $\theta_{c_{k}}=\vert \theta_{CT_{k}}-\theta_{TC_{k}}\vert$ is measured in radians.

\subsubsection{ Volume of ground cut}
If $s_{i-1},s_{i} \in T_{k,s}$ bracket a cut region, then
\begin{align}\label{CurveGcut}
a_{c}(s)=Wh_{c}(s)+\dfrac{1}{2}\kappa h_{c}^{2}(s),
\end{align}
where,
\begin{align}
h_{c}(s)= &\big( A_{uv}x_{c_{k}} +B_{uv}y_{c_k} +C_{uv} \big) + r_{k}\big(A_{uv}\cos\theta_{c_{k}}(s)+B_{uv}\sin\theta_{c_{k}}(s)\big) - \\ &\big( z_{k,{m_k}} + (z_{k+1,1} - z_{k,{m_k}})s   \big).
\end{align}
To calculate the Jacobian, $||r_c'(s)||$, we use
\begin{align}
r_c'(s) = \big(r_{k}\cos(\theta_{ck}(s))\theta_{ck}'(s),~r_k\sin(\theta_{ck}(s))\theta_{ck}'(s),~z_{(k+1),1} - z_{k,m_k} \big),
\end{align}
where $\theta_{ck}'(s) = \theta_{CT_k} - \theta_{TC_k}$. Hence, we have $$||r_c'(s)|| = \sqrt{r_k^2(\theta_{CT_k} - \theta_{TC_k})^2 + (z_{(k+1),1} - z_{k,m_k})^2}.$$

The volume of ground cut is, therefore, calculated as

\begin{align}
V_{c_k}^{i} = \left\{ \begin{array}{rl}
 	\textstyle \int_{s_{i-1}}^{s_{i}}a_{c}(s)||r_c'(s)||ds ~ & \mbox{if}~ s_{i-1}, s_i~ \mbox{bracket a cut segment}, \\
	0 & \mbox{otherwise} .
\end{array}\right.
\end{align}

The closed-form solution for $\int_{s_{i-1}}^{s_{i}}a_{c}(s)||r_c'(s)||ds $ can easily be calculated symbolically 
however, it requires almost full page to write. We therefore withhold it from this paper.

\subsubsection{Volume of ground fill}
If the parameters $s_{i-1}$ and $s_{i}$ bracket a fill region, the elevation difference at  $s\in [s_{i-1},s_{i}]$ is calculated as $h_{f}(s)=-h_{c}(s)$ and the cross-sectional area $a_{f}(s)$ is given by
 \[a_{f}(s)=Wh_{f}(s)+\dfrac{1}{2}\kappa h_{f}^{2}(s).\]
 Thus, the volume of ground fill is computed as

\begin{align}
V_{f_k}^{i} = \left\{ \begin{array}{rl}
 	\textstyle \int_{s_{i-1}}^{s_{i}}a_{f}(s)||r_c'(s)||ds ~ & \mbox{if}~ s_{i-1}, s_i~ \mbox{bracket a fill segment}, \\
	0 & \mbox{otherwise} .
\end{array}\right.
\end{align}


\subsection{Discussion on Simplifying Assumptions}

The model described above makes several simplifying assumptions with regards to the road design costs. We list some of these now, and mention that future research should explore the implications of these simplifying assumptions.  Up front, it should be clear that many road design costs are missing from the surrogate.  Earthwork cost is simplified to only include cut, fill, and borrow/waste costs.  Transportation costs are omitted, as are specific road design costs such as retaining walls and bridge work.  Utility cost is simplified even further, as it is represented using only the length of the road.  While many aspects of utility cost (such as paving, maintenance, speed-limits and travel time) are primarily affected by the length of the road, other factors can also affect these costs.  For example, maximum grade and road curvature can affect speed limits and travel time.  Like any real-world problem, these simplifications are made in order to create a computational tractable model.

Some more subtle simplifying assumptions are also present in the model.  For example, waste and borrow costs are combined into a single cost, $C_w$, which represents the amount of unbalanced material in the road design, $C_w \Vert V_{f}-V_{c}\Vert$.   In some situations one of the waste or borrow cost maybe significantly higher than the other.  In this case, the model could separate waste and borrow into two parameters by changing the first objective function to include $C_{b} \max\{0, V_{f}-V_{c}\} + C_{w} \max\{0, V_{f}-V_{c}\}$ instead of $C_{w}\Vert V_{f}-V_{c}\Vert$, where $C_b$ is the borrow cost and $C_w$ is the waste cost.

Another simplification made is the assumption that the earth shrinkage and earth swell factors are reciprocal. That is, overall no earth is lost or gained during the construction process.  In practice this is not always the case, and indeed will generally only be true if the optimized cut and fill quantities were exactly equal.  If one desired to incorporate shrinkage and swell factors into the model, it can be easily accomplished by multiplying the total cut and fill volume calculated above by an appropriate factor:
    $$V_c = \gamma_{\tt shrink} \sum_{\xi\in S}V_{c_{\xi}}$$
where $\gamma_{\tt shrink}$ is the shrinkage parameter.

Another simplification is the model assumes side slopes are equal at all points in the road.  Essentially this assumes that the ground material is of a similar composition at all point in the terrain.  To remove this assumption, the parameters $\theta_1$ and $\theta_2$ could be individually selected for each road segment.  Mathematically this will not change the model structure, or the computation of the integral, but in terms of implementation this would greatly complicate the code.

The model also omits transition curves from the road design, while past researchers have included transition curves \cite{Kang2012257}.  While similar approaches could be used to model transition curves, the model would require significant changes to accommodate them.

Also, when computing the cost of a road segment, the model computes the cost as if the center line of the road and the center line of the terrain is representative of the height of the road above (below) the terrain at the boundaries of the road.  If the terrain is sloped acutely with respect to the road, then it is possible that the terrain on one boundary of the road is above the road, while the terrain on the other boundary of the road is below the road. This would cause computational errors.  However, unless the road is particularly wide, the resulting error seems unlikely to be large.

Finally, the model ignores any technical design standards regarding interrelations between horizontal and vertical alignments.  This is a standard approach in road design literature, as it is assumed that an engineer will make final edits on any design.


\section{Bi-objective optimization model}\label{mop}
Multiobjective optimization (MOO) (also called multicriteria optimization,  or vector optimization) can be defined as the problem of finding a vector of decision variables which satisfies constraints and optimizes a vector of objective functions \cite{coello1999comprehensive}.  Bi-objective optimization refers to the specific case where the vector contains exactly two objective functions.

There is, usually, no unique solution that is simultaneously optimal for all objectives \cite{grana2005steepest}.  As a result, one can only consider a trade-off among the objectives, and the primary goal of multiobjective optimization is to seek for the best trade-off to support the decision maker  in choosing a final preferred solution. Although there is no universally accepted solution concept for decision problems with multiple  objectives, one would agree that a good solution must not be dominated by the other feasible alternatives \cite{yu1974cone}. The set of nondominated points is know as the \emph{Pareto optimal} set. Engineering design problems are often  multiobjective,  requiring trade-offs  \cite{wilson2001efficient}. In the literature, a great deal of theoretical, methodological, and applied studies have been undertaken in the area of multiobjective optimization \cite{figueira2005multiple}.

The need for MOO in road design has been recognized \cite{ATR5670430405,jha2007multi,Maji2011966}. During the road alignment design process, an engineer can have different objectives that need to be achieved. Some of the objectives may favour the shortest road possible, while others might favour an indirect and longer route with  smaller earthwork cost. The surrogate cost model developed in Section \ref{surrcost} has two components, the cost due to the volume of  earthwork and the cost related to the length of the road. There is, usually, a conflict between the two cost components.  We model the 3D road alignment optimization problem  as  a bi-objective optimization problem subject to the constraints listed in Section \ref{const}. The objectives are the utility cost $Cost_{\tt u}(X,Y,R,Z)$ and the earthwork cost $Cost_{\tt e}(X,Y,R,Z)$.   Thus, the solution of  the optimization problem should reflect the  trade-off between  the  length of the road and the volume of earthwork.

\subsection{Variable definition}
The costs are formulated as functions of decision variables. Let $X=\big(x_{1},x_{2},\cdots, x_{N}\big)$ and $ Y=\big(y_{1},y_{2},\cdots, y_{N}\big)$ be the coordinates of  the intersection points, $R=\big(r_{1},r_{2},\cdots, r_{N}\big)$ be the vector of  radius of curvature, and $Z=\big(z_{1},z_{2},\cdots, z_{M}\big)$ be the vector of elevations of the design road.

The input parameters of our model are: $START=(x_{s},y_{s},z_{s}),$  $END=(x_{e},y_{e},z_{e})$, and  the maximal vertical offset $\bar{z}$ from the current ground elevations $z_{g}$. The design parameters are the maximum allowable gradient $G_{\text{max}}$, and the minimum radius of curvature $r_{\min}$.
 Thus, the simultaneous optimization of horizontal and vertical alignments is given as follows.

\begin{align}\label{1}
  &\text{Minimize}~
   \begin{aligned}[t]
     \bigg\lbrace Cost_{\tt e}(X,Y,R,Z)~,~ Cost_{\tt u}(X,Y,R,Z)\bigg\rbrace
   \end{aligned} \notag\\
   &\text{subject to}:~\\
        &~~~~~~~~~~~~~~~~~\begin{aligned}[t]
        &\textbf{Horizontal alignment constraints}\\
        &\text{For}~ k=1,2,\cdots,N,\\
           &0\leq \Vert TC_{k}-IP_{k-1}\Vert -\Vert CT_{k-1}-IP_{k-1}\Vert, \\
&r_{\min}\leq r_{k}, \\
        &\textbf{Vertical alignment constraints}\\
        &\text{For}~ k=1,2,\cdots,N,~\text{ and for}~ j=1,2,\cdots m_{k},\\
        & \vert z_{k,j}-z_{k,j-1}\vert \leq d_{k,j}G_{\max},\\
        &z_{g}-\bar{z}\leq z_{k,j}\leq z_{g}+\bar{z}, \\
        &\textbf{Other constraints}\\
        &\text{For}~ k=1,2,\cdots,N,\\
        & x_{l,k}\leq x_{k}\leq x_{u,k},~  y_{l,k}\leq y_{k}\leq y_{u,k}
        \end{aligned} \notag
\end{align}


\section{Case Study}\label{numexp}

\begin{figure}[ht]
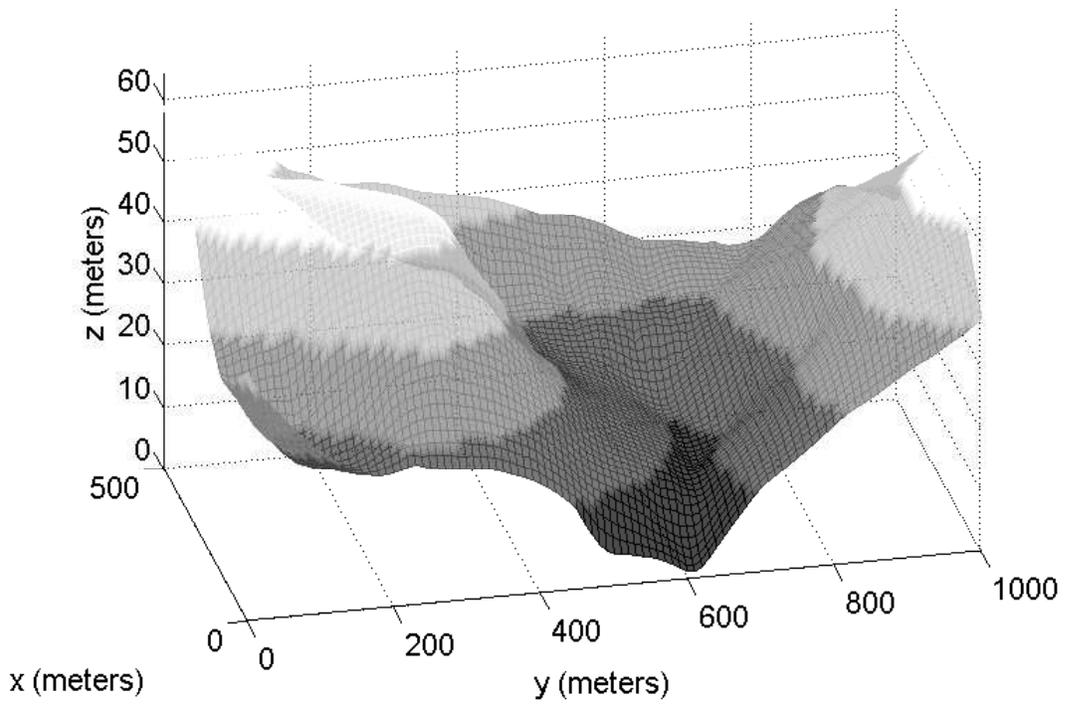
FIGURE 9 HERE\caption{3D view of the terrain}\label{fig30}\end{figure}

\begin{figure}[ht]
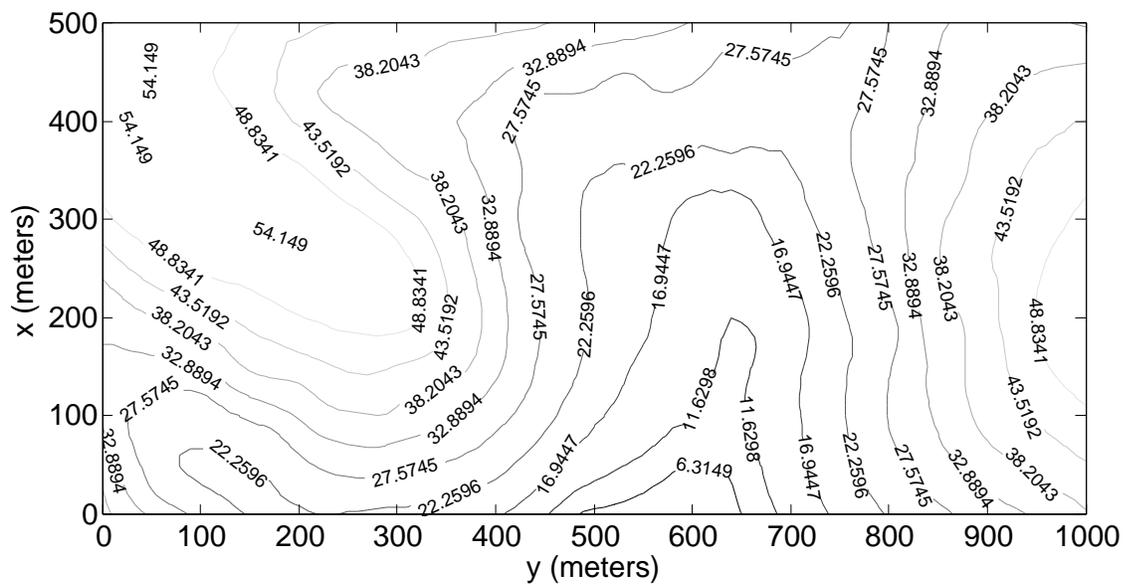
FIGURE 10 HERE\caption{Contour map of the terrain}\label{fig31}\end{figure}

In this section, we examine a case study applying the proposed model to design an approximately 1km section of road over an actual ground profile in California covering an area of 500 by 1000 meters. The 3D view and contour map of the terrain are displayed in Figure \ref{fig30} and Figure \ref{fig31}, respectively. As is typical in road design, the $z$ axis (road elevation) is measured in meters, while the $x,y$ axes (horizontal location) are measured in decameters. An initial alignment that satisfies all constraint sets is generated based on the ground profile of the study area.  The input values and design/model parameters are provided by our industrial partner, Softree Technical System Inc. (\url{ http://www.softree.com}) and based on approximate values for constructing a forest service road.  (Typical forest service roads cost between $\$10,000$ and $\$100,000$ per km.)  Cost details are presented in Table \ref{tab1}.

\begin{table}[ht]
	\centering
	\caption{Input variables and design/model parameters.}\label{tab1}
	\begin{tabular}{ c c c c c c c c c c c}
		\hline
		Parameter \vline & C$_{u}$ & C$_{c}$ &  C$_{f}$ & C$_{w}$ & G$_{\max}$ (\%)  & r$_{\min}$ (m) & $m_{k}$ & N & W & $\kappa$    \\ [0.5ex]
		\hline
		Values ~~~~~\vline & 1.2 & 4 & 2 & 8 & 15 & 20 & 5& 6 & 5 & 1\\ [0.5ex]
		\hline
	\end{tabular}
\end{table}

Numerical experiments were designed and conducted with a MATLAB R2013b code performed on a Dell workstation equipped with an Intel(R) Xeon(R) CPU E5-1620 v2 3.70GHz processor, and 32 GB of RAM using 64-bit Windows operating system.  A test problem is solved  using  three different  optimization algorithms namely, the  multiobjective genetic algorithm (MOGA) in MATLAB's Global Optimization Toolbox, the direct multisearch for multiobjective optimization (DMS) \cite{dmsMOO}, and the weighted sum method (WS) \cite{normalization, weightedsum}.

\subsection{Multiobjective genetic algorithm}
A genetic algorithm (GA) is a search algorithm inspired by the principle of natural selection.
The basic idea is to evolve a population of individuals, which are evaluated by a fitness function that measures the quality of its corresponding solution. At each generation (iteration) the fittest (the best) individuals of the current
population survive and produce offspring resembling them, so that the population gradually contains better individuals. The GAs have found application  in road design optimization problems   which are known to present difficulties to  conventional numerical optimization \cite{MICEMICE574, MICEMICE626, MICEMICE778, jha2007multi,Jong2003107}.

A  MOGA solver in MATLAB software, which is based on  the concept of Pareto dominance, is employed to solve the example scenario.

\begin{figure}[ht]
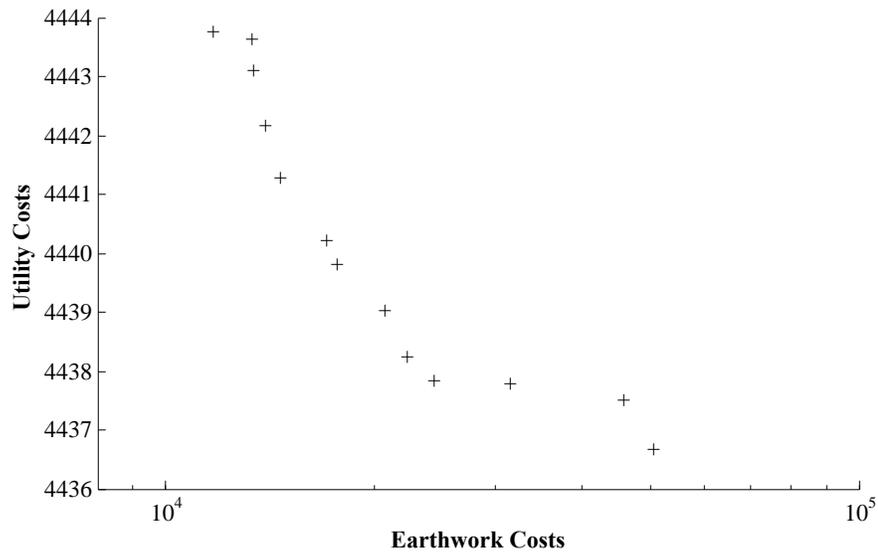
FIGURE 11 HERE  \caption{Pareto front using MOGA with \emph{TolFun}=$10^{-4}$.} \label{fig:fig34}\end{figure}

Based on the recommendation given in \cite{deb2001multi}, we set the population size to 120 and set \emph{TolFun} to $10^{-4}$.
The resulting Pareto front is depicted in Figure \ref{fig:fig34}.

\subsection{Direct multisearch for  multiobjective optimization}

The direct multisearch (DMS) for  multiobjective optimization is a derivative-free optimization algorithm that uses the concept of Pareto dominance to maintain a list of nondominated points. The method does not aggregate any of the objective functions, instead it  extends the classical directional derivative-fee methods from single to multiobjective optimization problem \cite{dmsMOO}.

\begin{figure}[ht]
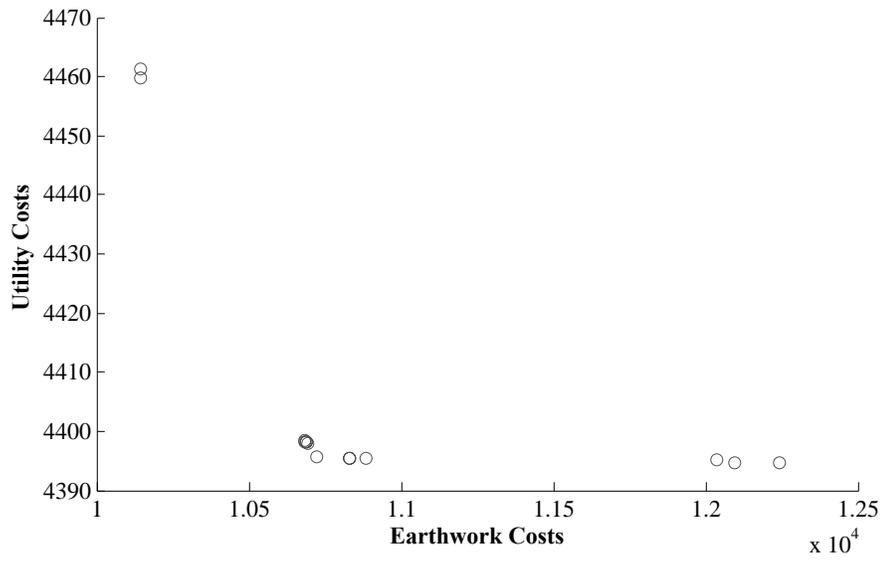
FIGURE 12 HERE\caption{Pareto front using DMS with\\ $f$-count=51,001}\label{fig36}\end{figure}

In order to compare the results of the three solvers, we use the number of  function calls as a stopping criteria for the other solvers. The amount of function calls used by MOGA was 51,231, hence we run each remaining solver until 51,000 function calls have been surpassed (note that function calls are checked at the end of an iteration, so as a result solvers will use over 51,000 function calls).  DMS required 51,001 function calls to generate its Pareto front, which is shown in Figure \ref{fig36}.

In comparing figures \ref{fig:fig34} and \ref{fig36}, note that scales differ.  A clean figure, comparing all three solvers, appears in Subsection \ref{ss:compare}.

\subsection{Solution by the weighted sum method }

The weighted sum (WS) method scalarizes a set of objectives into a single objective by  pre-assigning each objective with a user-supplied weight \cite{weightedsum}. As different objective functions can have different magnitude, the normalization of objectives is required to get a Pareto optimal solution consistent with the weights assigned by the user \cite{normalization}. For the purpose of  normalization, we use a \emph{Nadir} point $C^{N}\in \mathbb{R}^{2}$ and a \emph{Utopia} (or \emph{Ideal}) point $C^{I}\in \mathbb{R}^{2}$. A Nadir point is  defined as the vector whose components are the individual maxima in the Pareto front of each objective function, and a Utopia point is  the vector whose components are the individual minima of each objective function \cite{audet2008multiobjective}.
So, the weights are computed as

\begin{equation}\label{singleobject}
w_{e}=v_{e} N_{e}, ~~w_{u}=v_{u} N_{u},
\end{equation}
 where $0\leq v_{e}\leq 1, v_{u}=1-v_{e}$ are user-supplied weights, and  $N_{e}$, $N_{u}$ are the normalization factors  calculated as follows \cite{normalization}.
\[N_{e}=\dfrac{1}{C_{e}^{N}-C_{e}^{I}},~ N_{u}=\dfrac{1}{C_{u}^{N}-C_{u}^{I}},\]
where $C_{u}^{N}, C_{e}^{N}$ are components of the Nadir point and $C_{e}^{I}, C_{u}^{I} $ are components of the Ideal point.

We solved the scalarized objective functions defined by
\[\text{Cost}~=w_{e}Cost_{\text{e}}+w_{u}Cost_{\text{u}}\]
for 51 values of the weight $v_{e}$ (running between 0 and 1 with  a step of $1/50$) using the  solver \emph{fmincon} in  MATLAB. The  \emph{interior-point}  algorithm is used to solve the problem. The  number of function calls was used as a stopping criteria.  In this case, we run only one experiment by setting  the number of  function counts to 1,000.  The values of objective functions at each  solution point  are plotted, see Figure \ref{fig37}, where, the boxed points are nondominated points that we identified.

\begin{figure}[ht]FIGURE 13 HERE\caption{All points found using the Weighted sum method. Boxed points representing the Pareto front.}\label{fig37}
\end{figure}

The number of non-dominated points that we identified from the weighted sum method are few in number. This is because the scalarized objective function is solved only 51 times and many of the solutions are dominated.   We did not set any  special procedure to sort a list of non-dominated points. Instead, we evaluated each objective function at a solution point that we obtain by solving the scalarized objective, and the   objective functions values are plotted in the objective space. Then, the set of non-dominated points are identified  by comparing the objective function values.

\subsection{Comparison of Methods and Engineering Implications}\label{ss:compare}
\begin{figure}[ht]
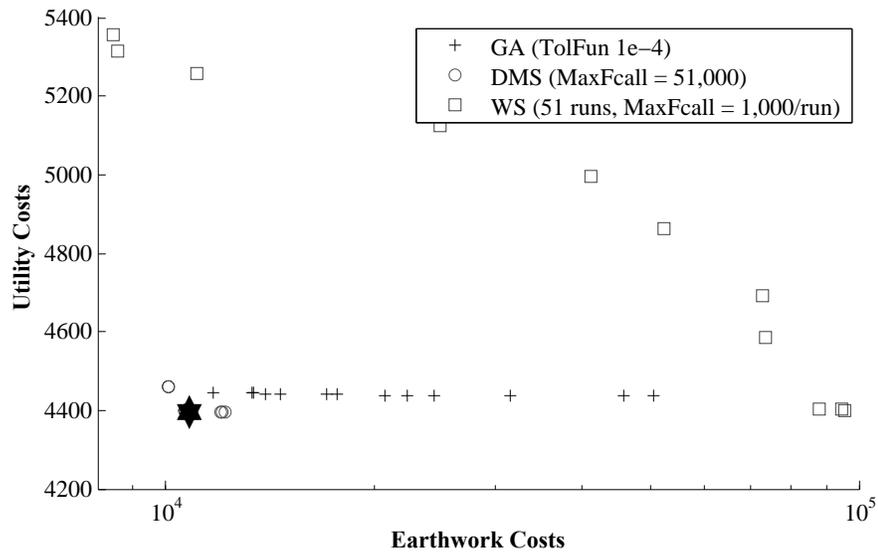
FIGURE 14 HERE\caption{Pareto fronts from MOGA, DMS, and WS, plotted on same axes. The black star is an example solution seen in Figure~\ref{fig:SampleRoad}.}\label{fig38}\end{figure}

In Table \ref{tab3}, we summarized the results of numerical experiments. The Pareto front of all test cases is plotted in Figure \ref{fig38}.

\begin{table}
\centering
\caption{Summary of numerical experiments.}\label{tab3}
\begin{tabular}{ c c c c}
 \hline
 Solver &   f-counts  & computation time (sec.)  \\  [0.5ex]
 \hline
 GA (\emph{TolFun}=$1e^{-4}$) & ~~~~51,231  & 729  \\
 DMS (\emph{MaxFcall}=$51,000$) & ~~~~51,001 & 795 \\
 WS (51 runs, \emph{MaxFcall}=$1,000$/run)  & ~~~~52,858  & 802 \\
 \hline
\end{tabular}
\end{table}

\begin{figure}[ht]
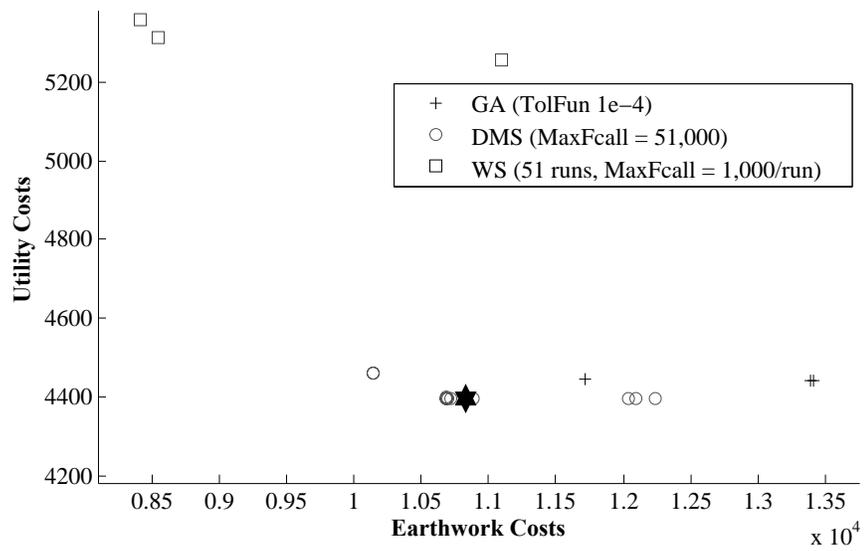
FIGURE 15 HERE\caption{The zoomed in Pareto fronts from MOGA, DMS, and WS in Figure~\ref{fig38}. The black star is an example solution seen in Figure~\ref{fig:SampleRoad}.}\label{fig:Zoomed All in one}\end{figure}

In examining Figure \ref{fig38}, note that the earthwork cost ranges from $10^4$ to $10^5$, while the utility cost only ranges from $4200$ to $5400$. This should not be taken to mean that earthwork costs are more than utility costs.  As mentioned before, earthwork costs can be viewed as present day dollar (the `construction dollars'); but utility costs should be viewed as a value representing long term road effectiveness, so does not necessarily have a clear translation to present day dollars.  What is clear in Figure \ref{fig38} is that the high end of the earthwork costs shows very little improvement in utility costs.  As such, in Figure~\ref{fig:Zoomed All in one} we present a cropped view of the Pareto front, focusing on earthwork costs of $8000$ to $13750$.

The `starred' point in Figures~\ref{fig38} and \ref{fig:Zoomed All in one} was found using DMS and corresponds to a utility cost of $4395$ and an earthwork cost of $10829$. Based on the numerical results obtained, we see the DMS performs better than others in terms of the magnitude of objective functions when the function call is set to about 51,000.  While MOGA provides a large spread in the Pareto front, all of the points found are Pareto dominated by the starred point.

\begin{figure}[ht]
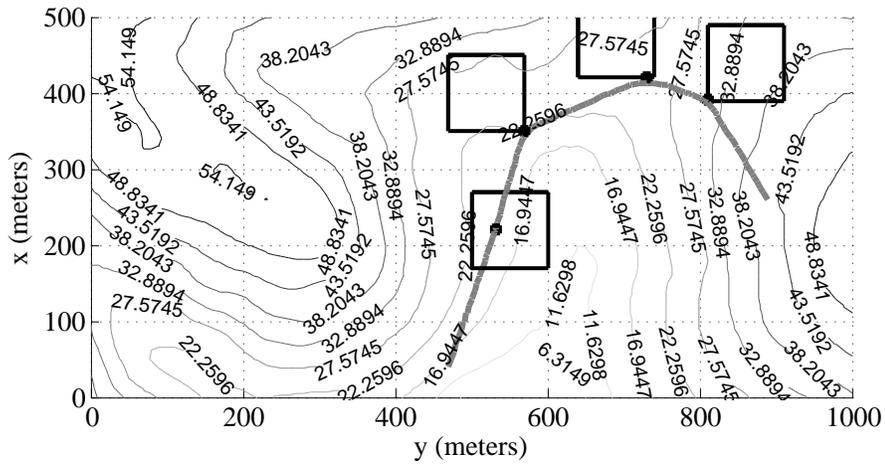
FIGURE 16 HERE\caption{The horizontal alignment corresponding to the starred Pareto optimal point in Figure \ref{fig38}. The black squares represent the box constraints on the intersection points and the block points represent the intersection points.} \label{fig:SampleRoad}\end{figure}

Working from Figure~\ref{fig:Zoomed All in one}, we conclude that, in this case study, there are a relatively small number of roads that require further investigation when selecting a final road design.  In Figure \ref{fig:SampleRoad}, we present the horizontal alignment for the road corresponding to the starred point in Figures \ref{fig38} and \ref{fig:Zoomed All in one}.

\section{Conclusion}\label{colcln}
 In this study, a  model  was developed  to solve a 3-dimensional road alignment optimization problem. The model uses bi-objective optimization to minimize the earthwork and utility costs. Cost penalty parameters are introduced and their values (that we obtain from our industry partner) are fixed, but, in theory, they can be computed in order to calibrate  the solution of the  model. The cost items are classified into those that depend on the length of the road and those that depend on the volume of the earthwork. These are, often, conflicting with each other because, during the optimization process, the utility cost prefers the shortest route between the two end points, while the earthwork cost chooses a route with minimum amount of earthwork.

 As a major contribution to the field of three-dimensional alignment, the paper provides a high level of detail on exactly how the cost computations are preformed.  The model makes a number of simplifying assumptions, discussed in detail in Section \ref{}.  It is our hope that the detailed cost computations outlined will allow for improved accessiblity for future research to explore the impact of these assumptions and remove them when deemed necessary.

 A case study based on actual terrain was examined. Three different optimization algorithms (MOGA, DMS, and WS) were used to solve the bi-objective optimization problem. In this case study, DMS provided the best Pareto front for this model. More importantly, both MOGA and DMS were capable of providing a good quality Pareto front within a reasonable time limit.  Further research is needed to check whether MOGA and DMS are consistent in this regard.

 It should be noted that MOGA, DMS, and WS were chosen due to availability of software, and therefore are not necessarily the best solvers for this problem.  Other researchers have explored creating custom solvers for road design optimization (\cite{schonfeld2006intelligent, Jong2003107, ATR5670430405} among others).  Future work will progress in this direction.

\section*{Acknowledgments}

This work was supported by the Natural Sciences and Engineering Research Council of Canada (NSERC) through Collaborative Research and Development grant \#CRDPJ 411318-2010 sponsored by Softree Technical Systems Ltd. and Discovery
grants \#355571-2013 (Hare), \#298145-2013 (Lucet), and \#2014-05013 (Tesfamariam). Part of the research was performed in the Computer-Aided Convex Analysis (CA2) laboratory funded by a Leaders Opportunity Fund (LOF) from the Canada Foundation for Innovation (CFI) and by a British Columbia Knowledge Development Fund (BCKDF).  Hirpa received support from the University of British Columbia (UBC) through University Graduate Fellowships (UGF).

The authors would like to acknowledge the assistance and positive impact of the anonymous referees who helped improve this work.


~\newpage~
\setcounter{figure}{0}

\begin{figure}
\centering
\includegraphics[scale=0.3]{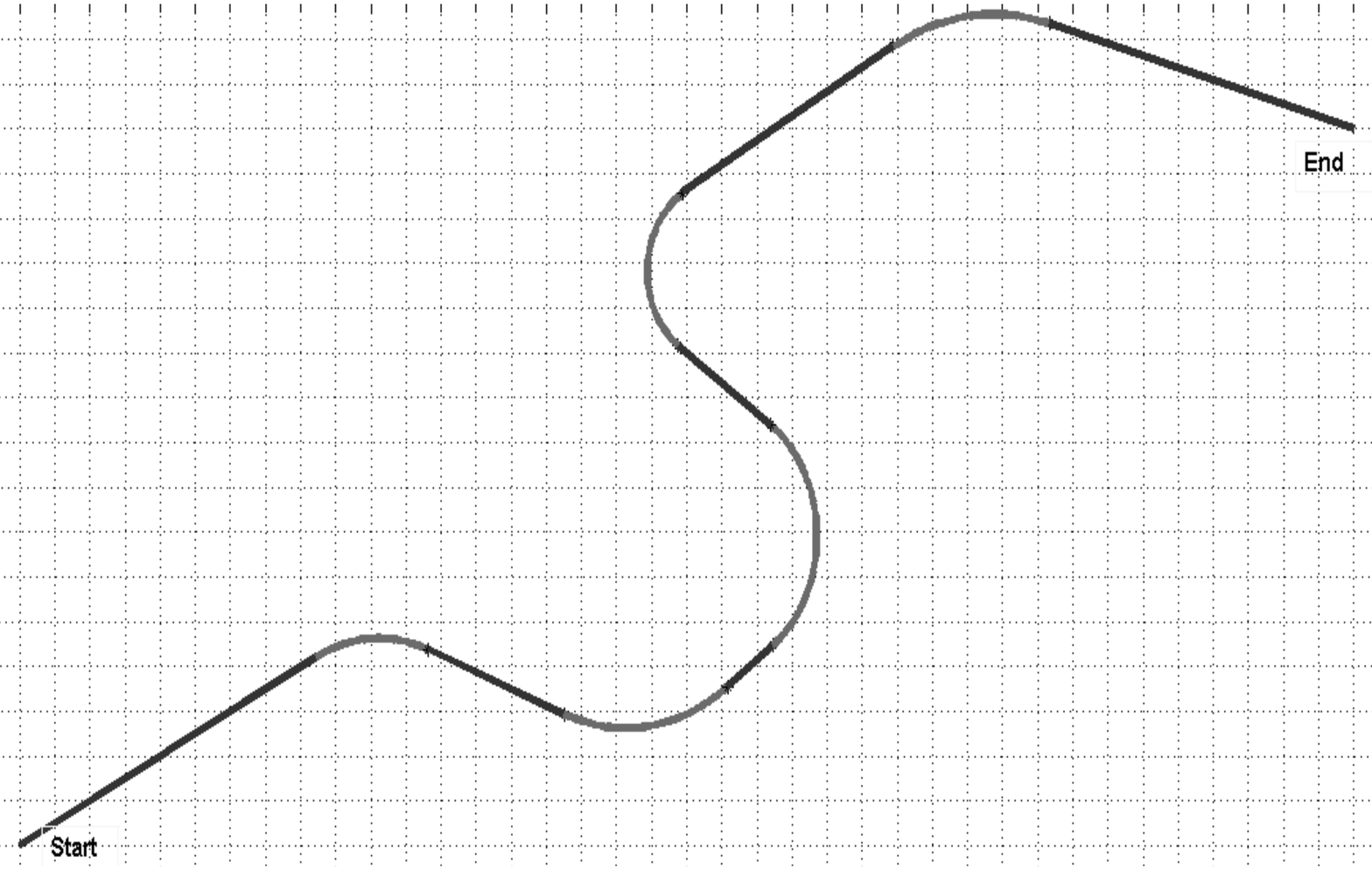}
\caption{An example of  horizontal alignment.}
\end{figure}

~~\newpage~~

\begin{figure}
\centering
\includegraphics[scale=0.4]{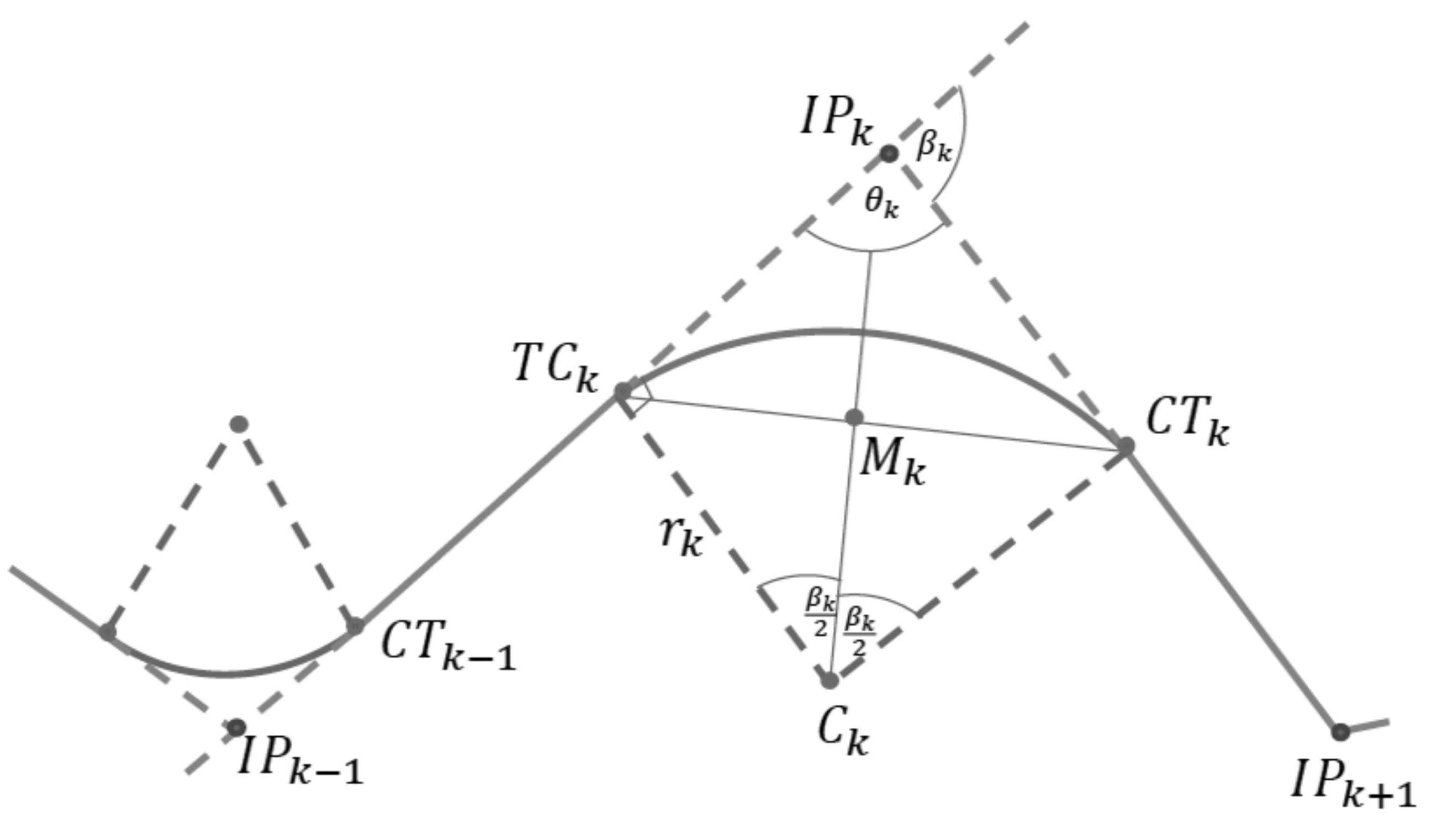}
\caption{Horizontal alignment geometry}
\end{figure}

~\newpage~

\begin{figure}
\centering
\includegraphics[scale=0.4]{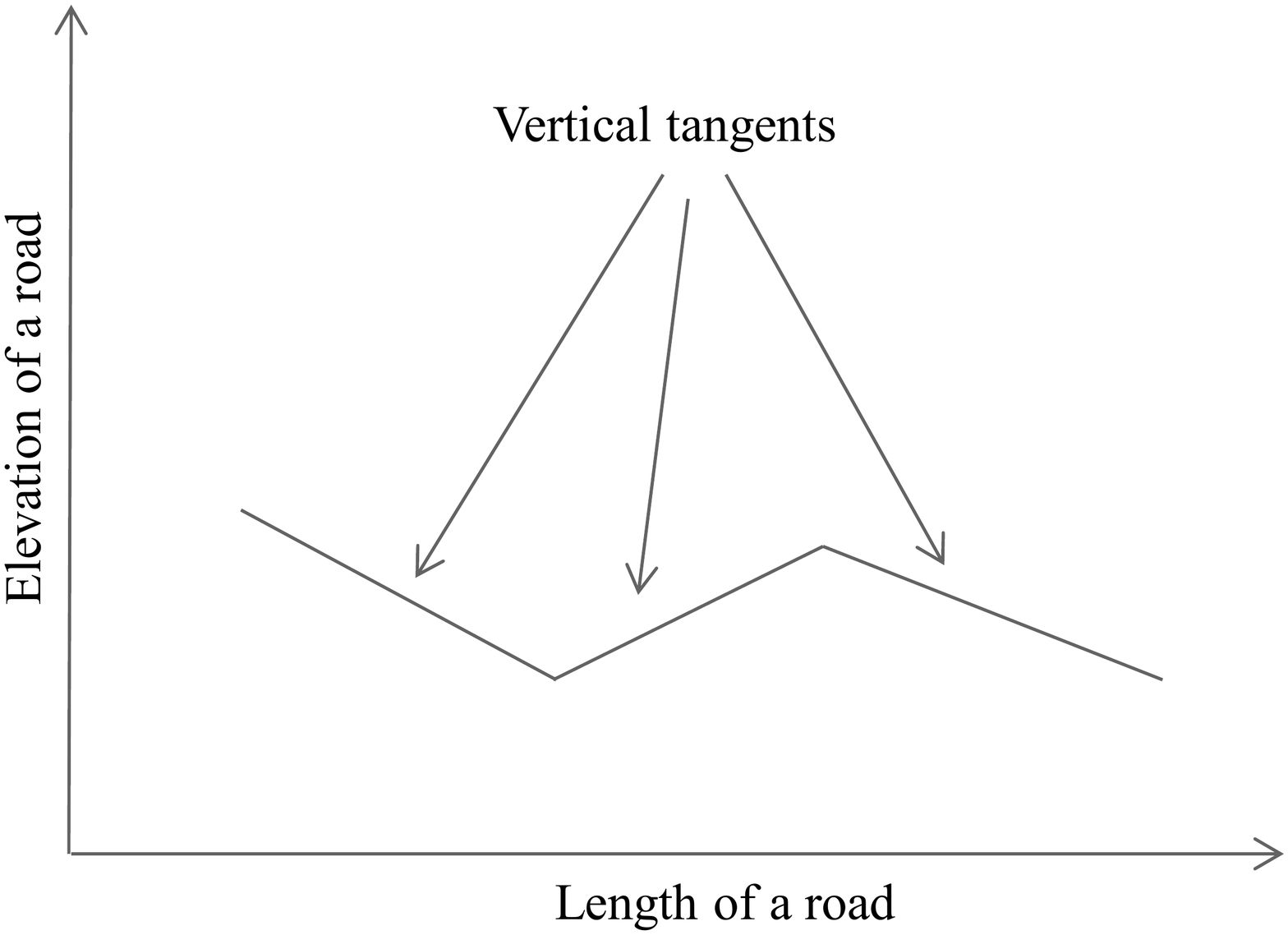}
\caption{An example of   vertical alignment model.}
\end{figure}

~\newpage~

\begin{figure}
 \begin{center}
\includegraphics[scale=0.4]{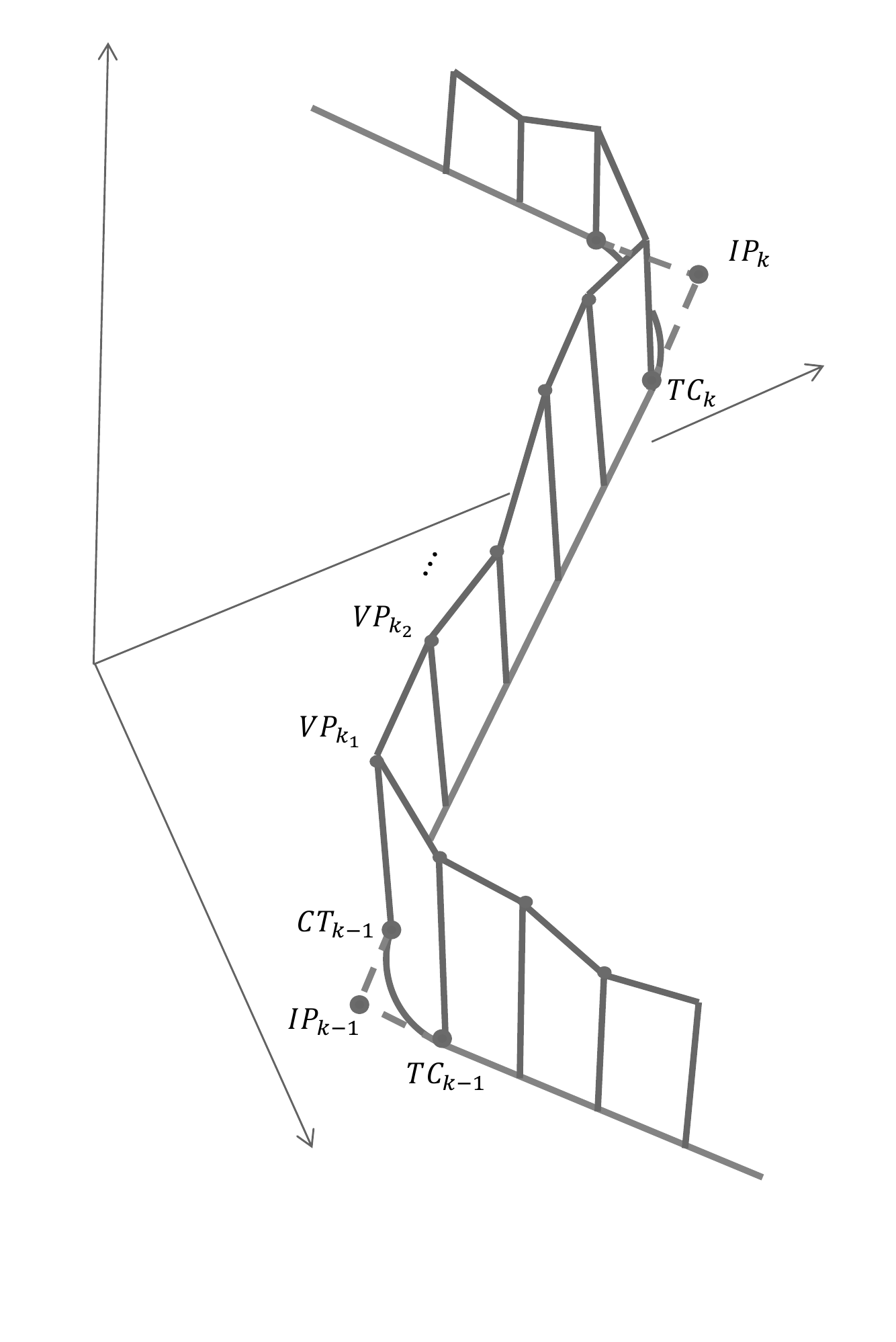}
\caption{Road alignment in 3D.}
\end{center}
\end{figure}

~\newpage~

\begin{figure}
\centering
\includegraphics[scale=0.35]{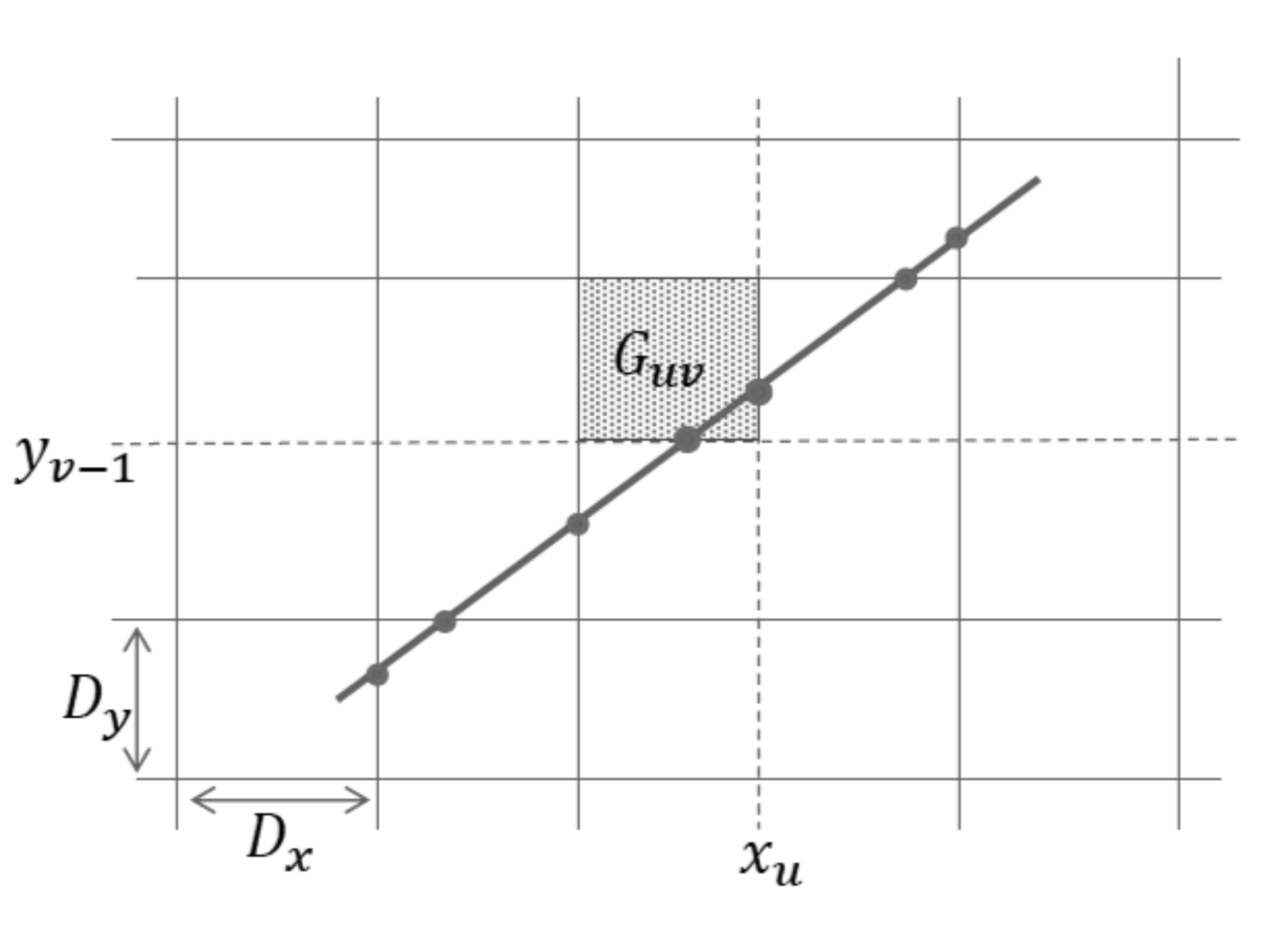}
\caption{An example, the projection of a tangent road segment onto the horizontal plane.}
\end{figure}

~\newpage~

\begin{figure}
\centering
\includegraphics[scale=0.4]{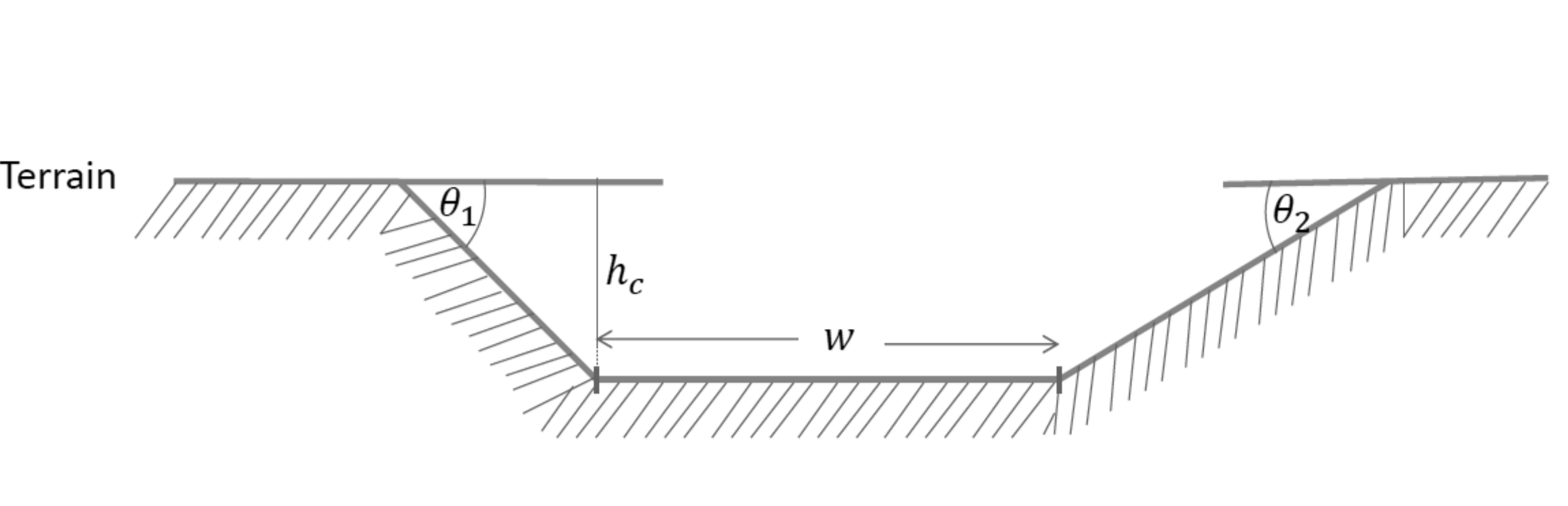}
\caption{An example of cut cross-section. }
\end{figure}

~\newpage~

\begin{figure}
\centering
\includegraphics[scale=0.4]{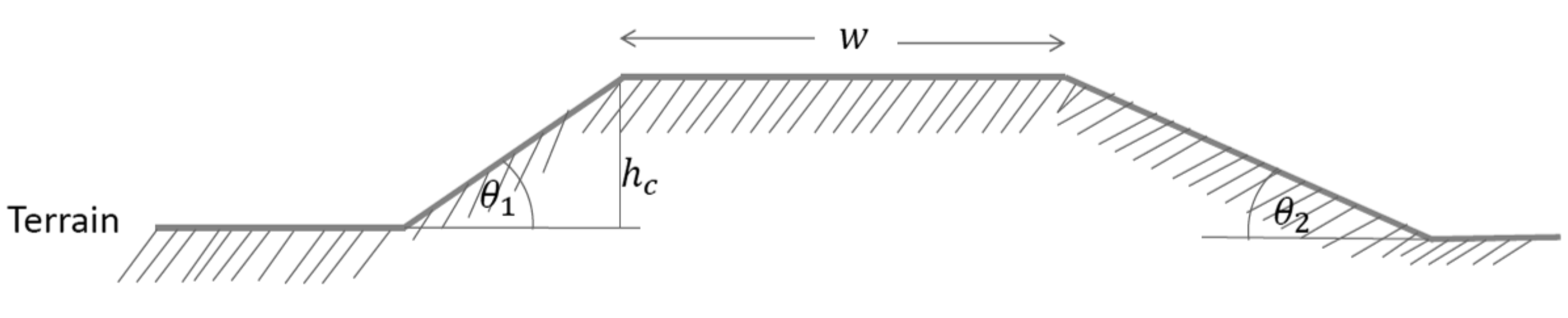}
\caption{An example, fill cross-section. }
\end{figure}

~\newpage~

\begin{figure}
\centering
\includegraphics[scale=0.5]{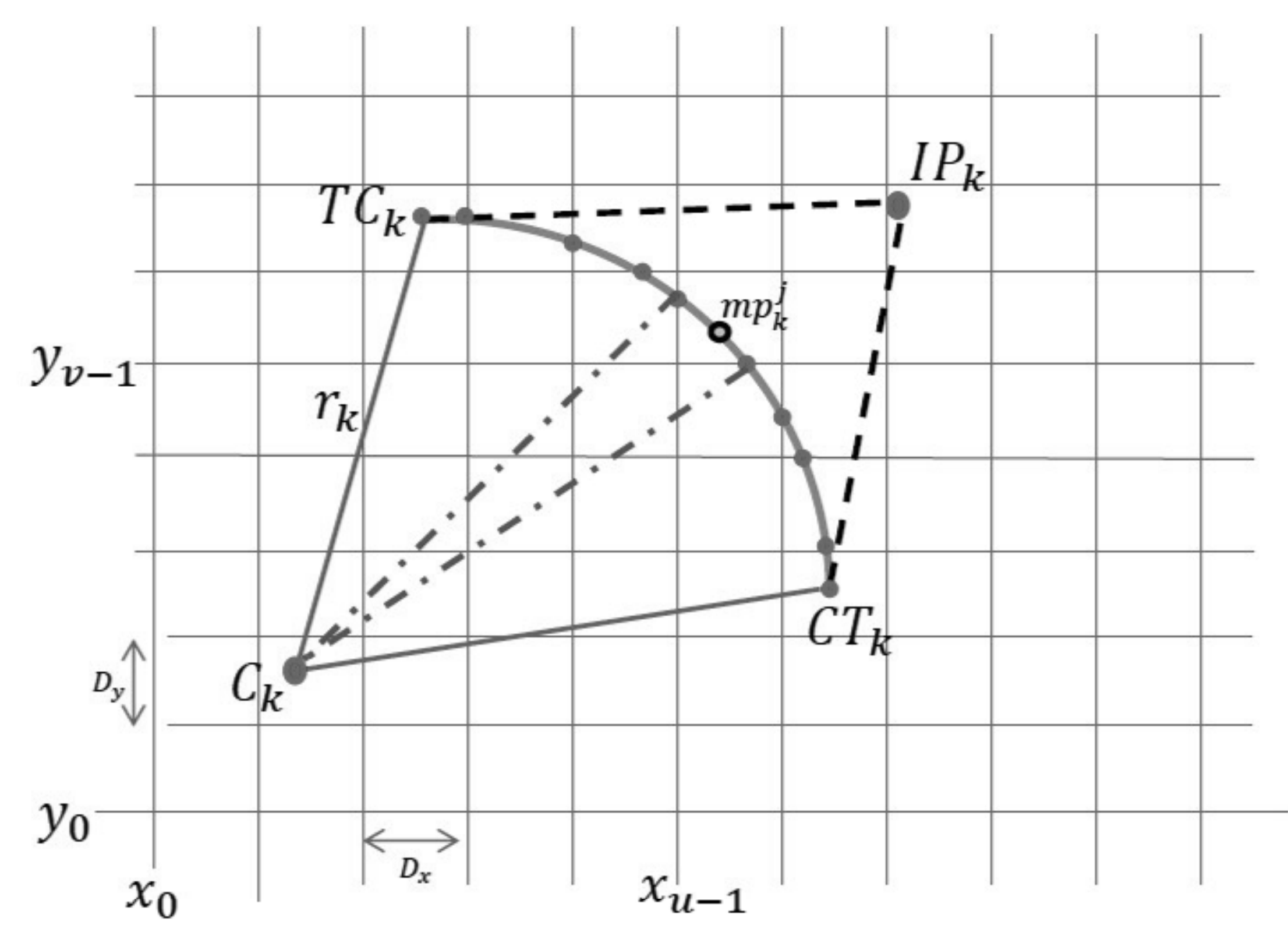}
\caption{Example of a horizontal curve section}
\end{figure}

~\newpage~

\begin{figure}[ht]
\centering
\includegraphics[scale=0.7]{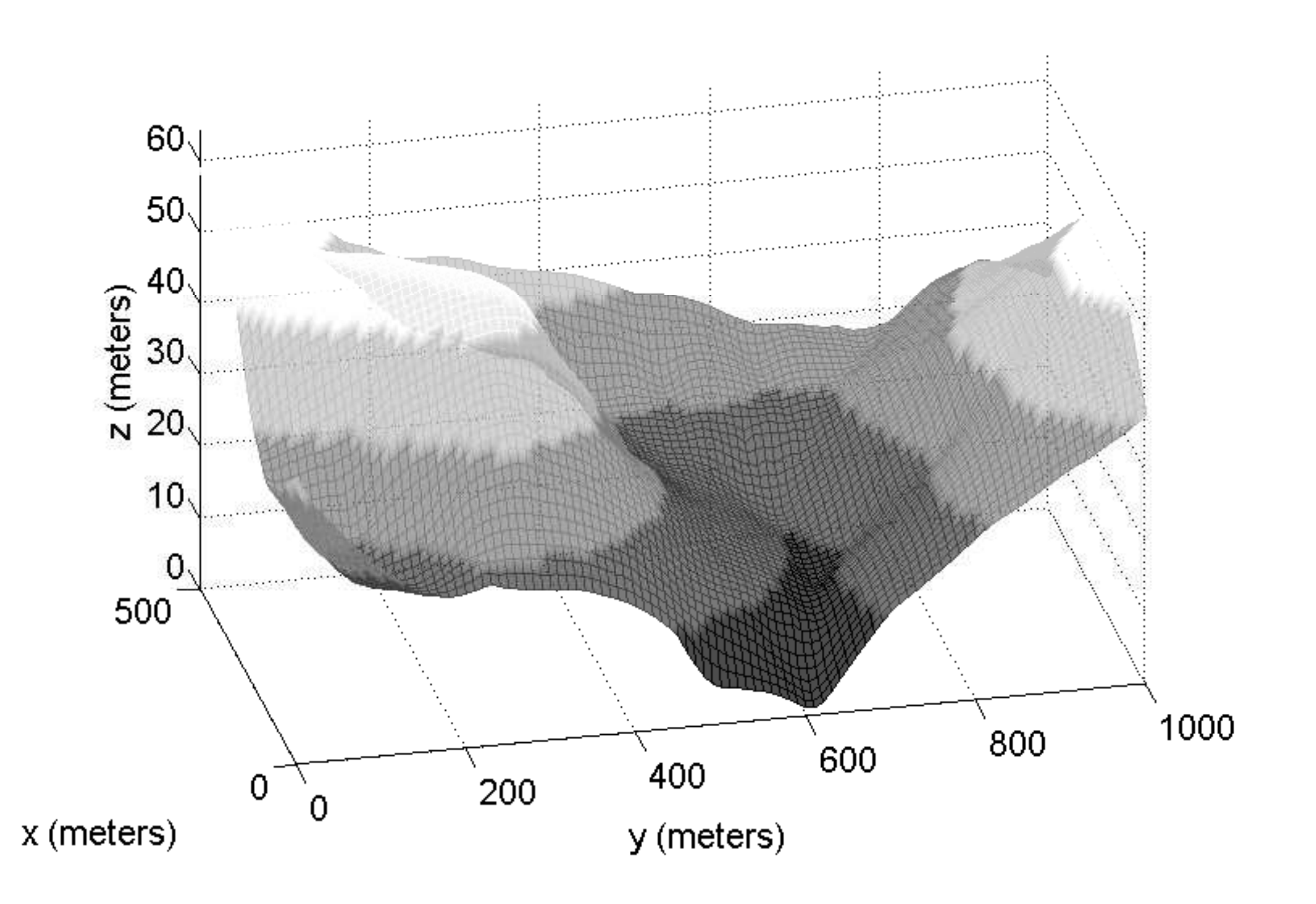}
\caption{3D view of the terrain}
\end{figure}

~\newpage~

\begin{figure}[ht]
\centering
\includegraphics[scale=0.7]{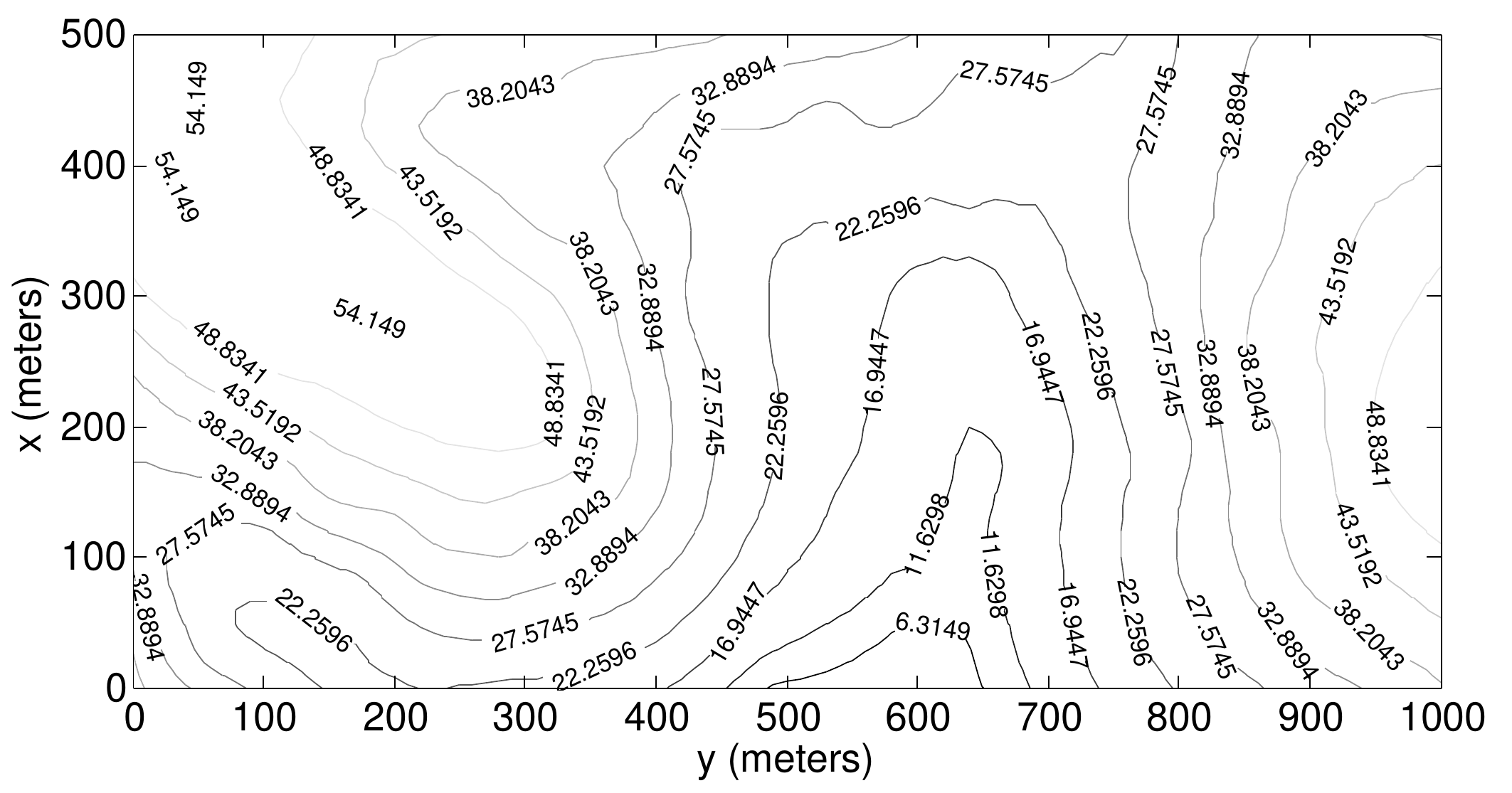}
\caption{Contour map of the terrain}
\end{figure}

~\newpage~

\begin{figure}
  \centering
  \includegraphics[scale=0.75]{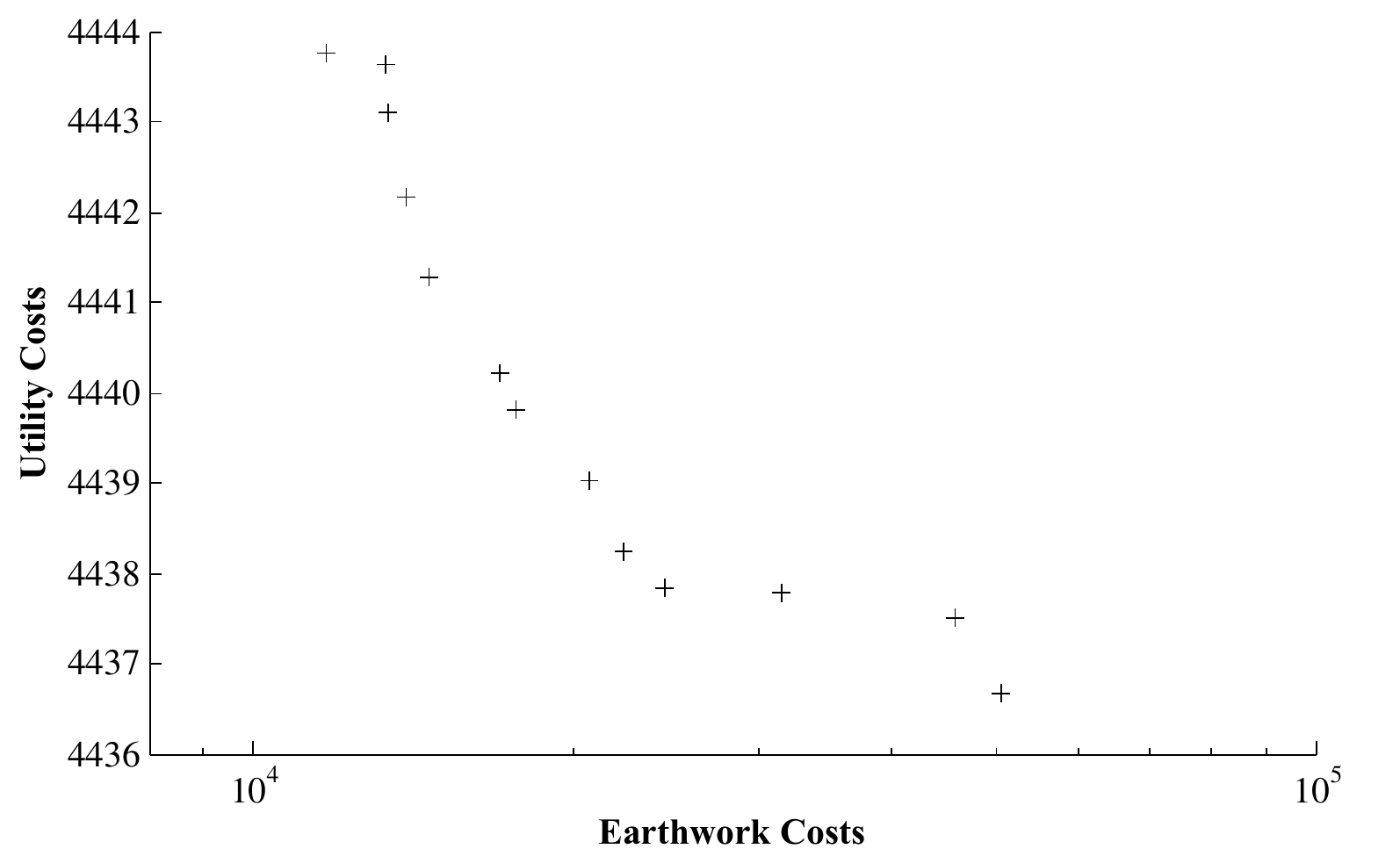}
  \caption{Pareto front using MOGA with \emph{TolFun}=$10^{-4}$.}
\end{figure}

~\newpage~

\begin{figure}
	\centering
	\includegraphics[scale=0.75]{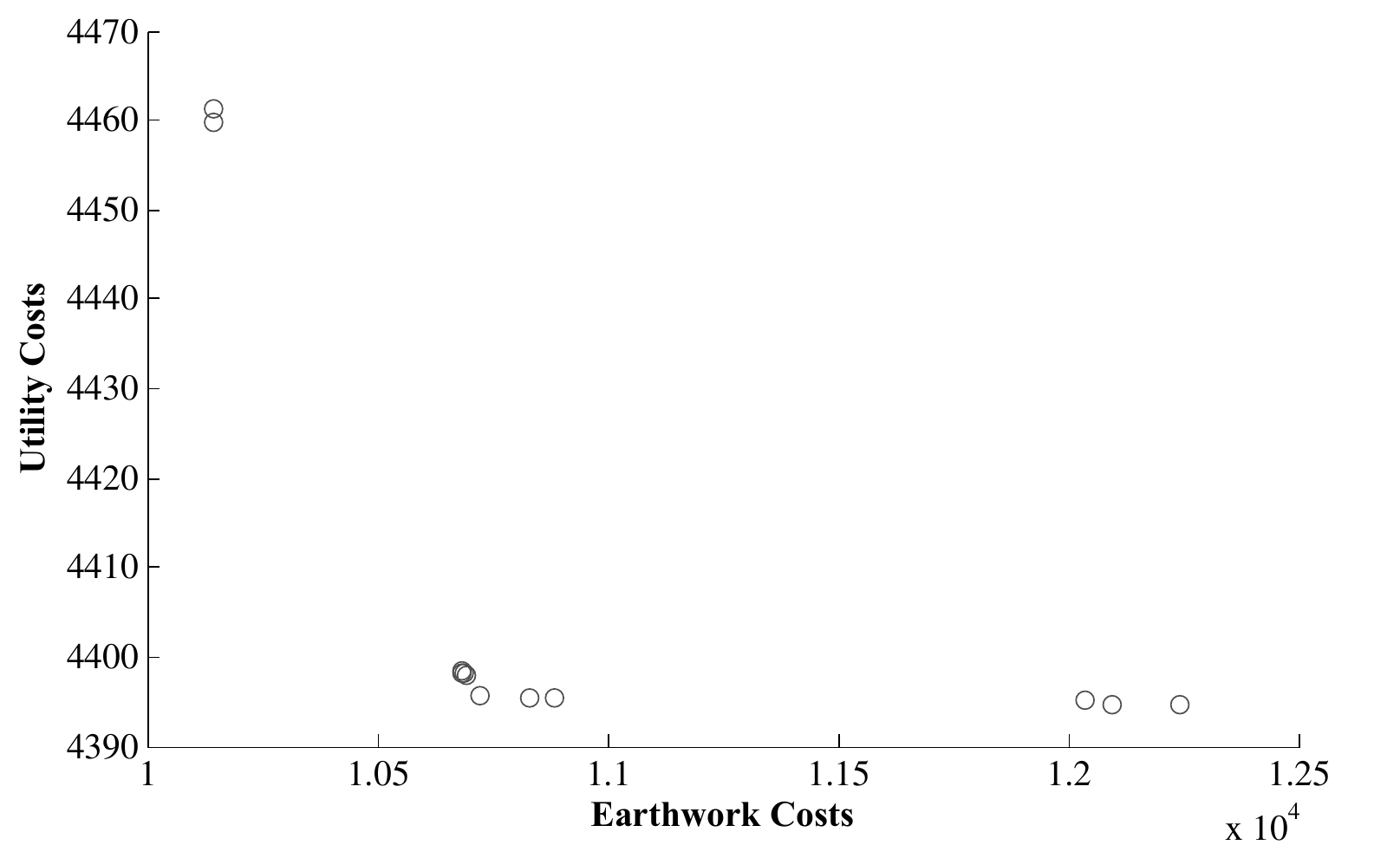}
	\caption{Pareto front using DMS with\\ $f$-count=51,001}
\end{figure}

~\newpage~

\begin{figure}
\centering
\includegraphics[scale=0.75]{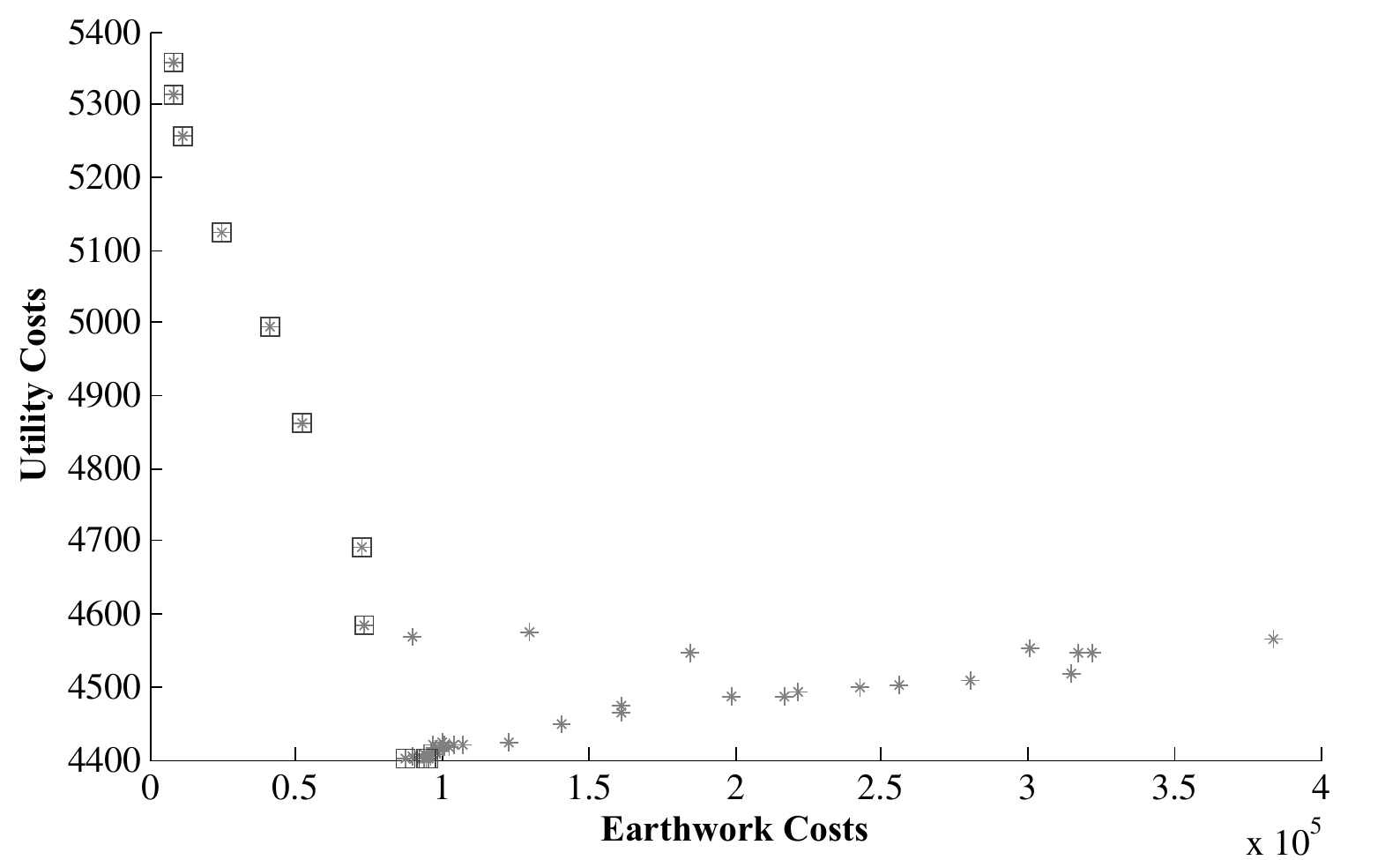}
\caption{All points found using the Weighted sum method. Boxed points representing the Pareto front.}
\end{figure}

~\newpage~

\begin{figure}
\centering
\includegraphics[scale=0.75]{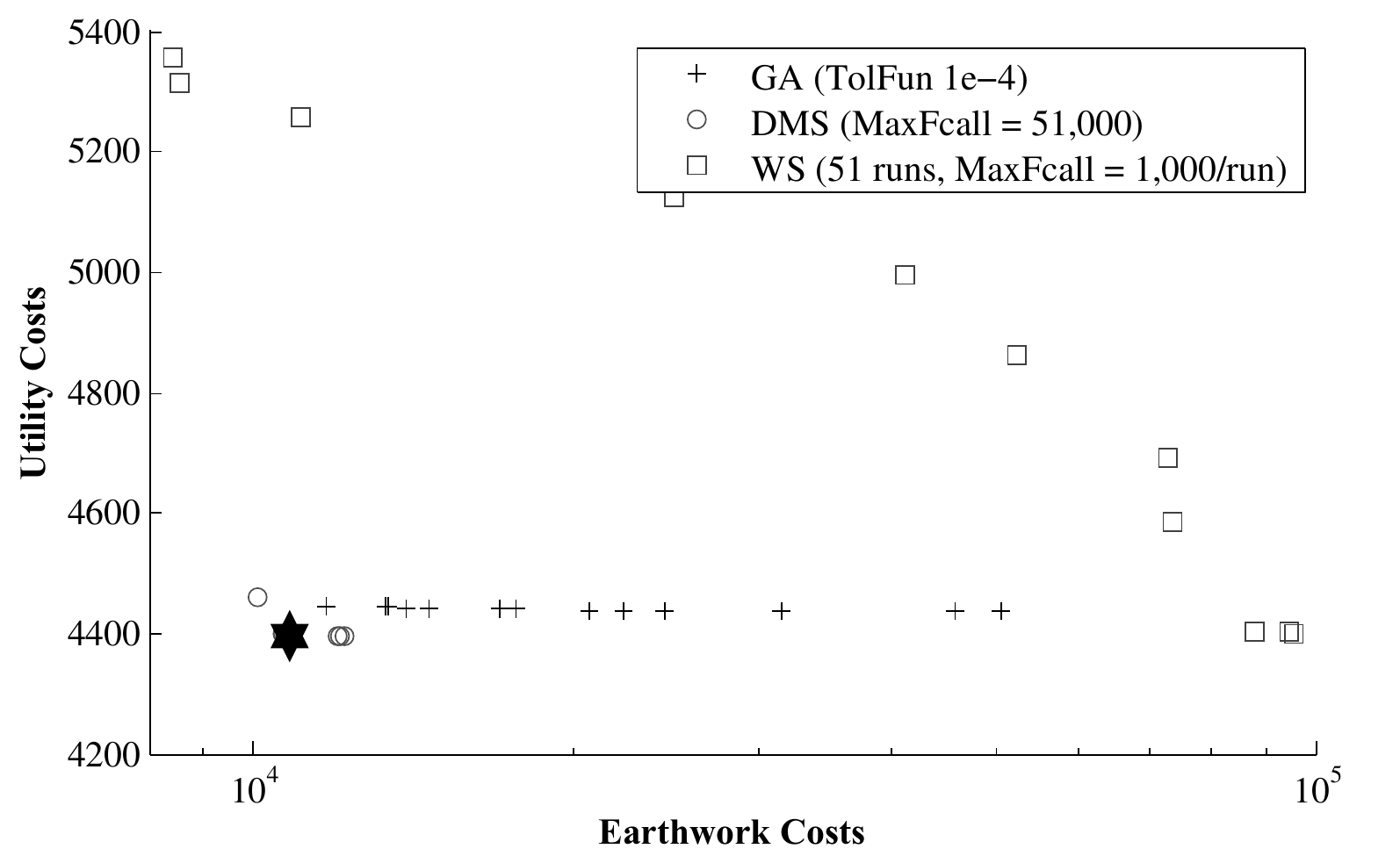}
\caption{Pareto fronts from MOGA, DMS, and WS, plotted on same axes. The black star is an example solution seen in Figure~\ref{fig:SampleRoad}.}
\end{figure}

~\newpage~

\begin{figure}
	\centering
	\includegraphics[scale=0.75]{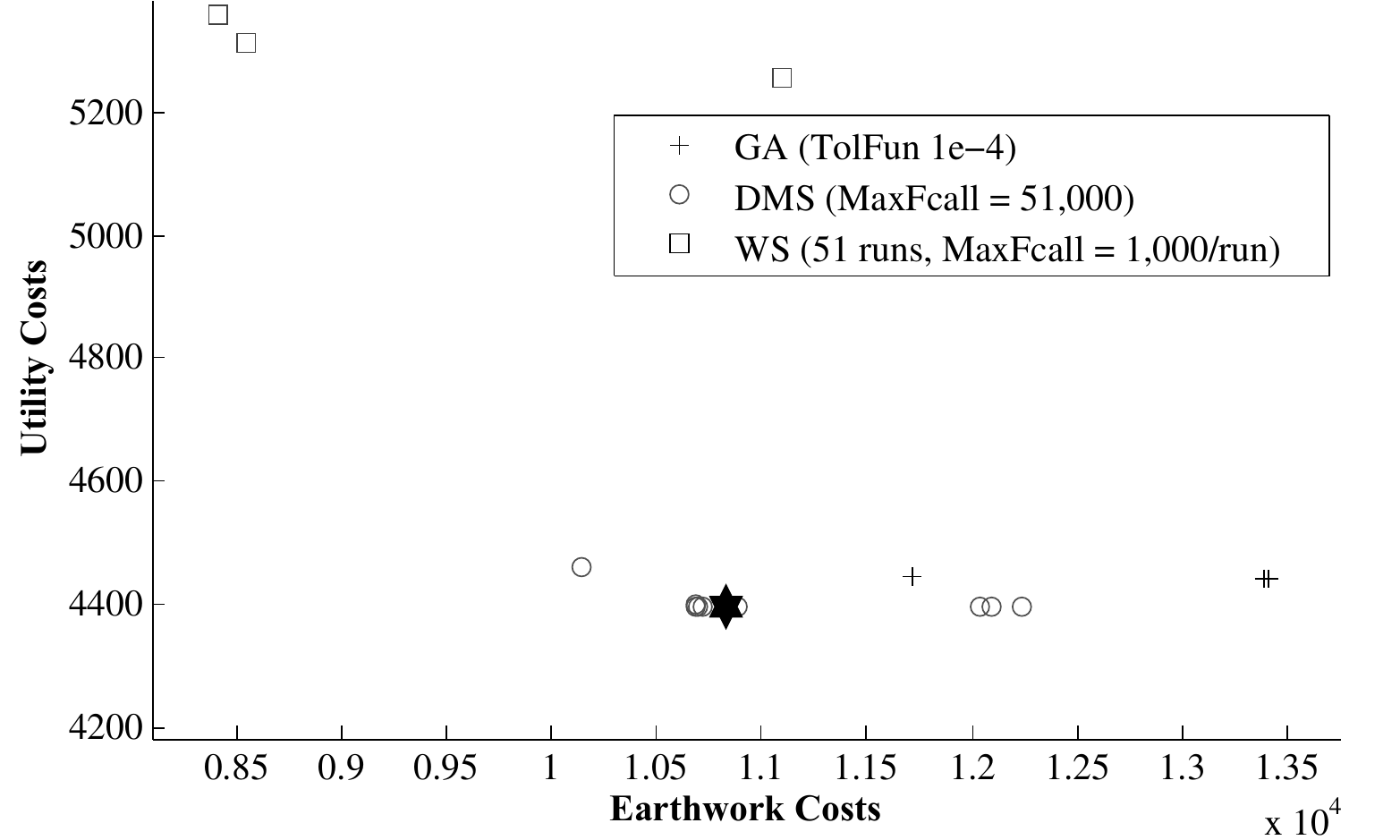}
	\caption{The zoomed in Pareto fronts from MOGA, DMS, and WS in Figure~\ref{fig38}. The black star is an example solution seen in Figure~\ref{fig:SampleRoad}.}
\end{figure}

~\newpage~

\begin{figure}
	\centering
	\includegraphics[scale=0.75]{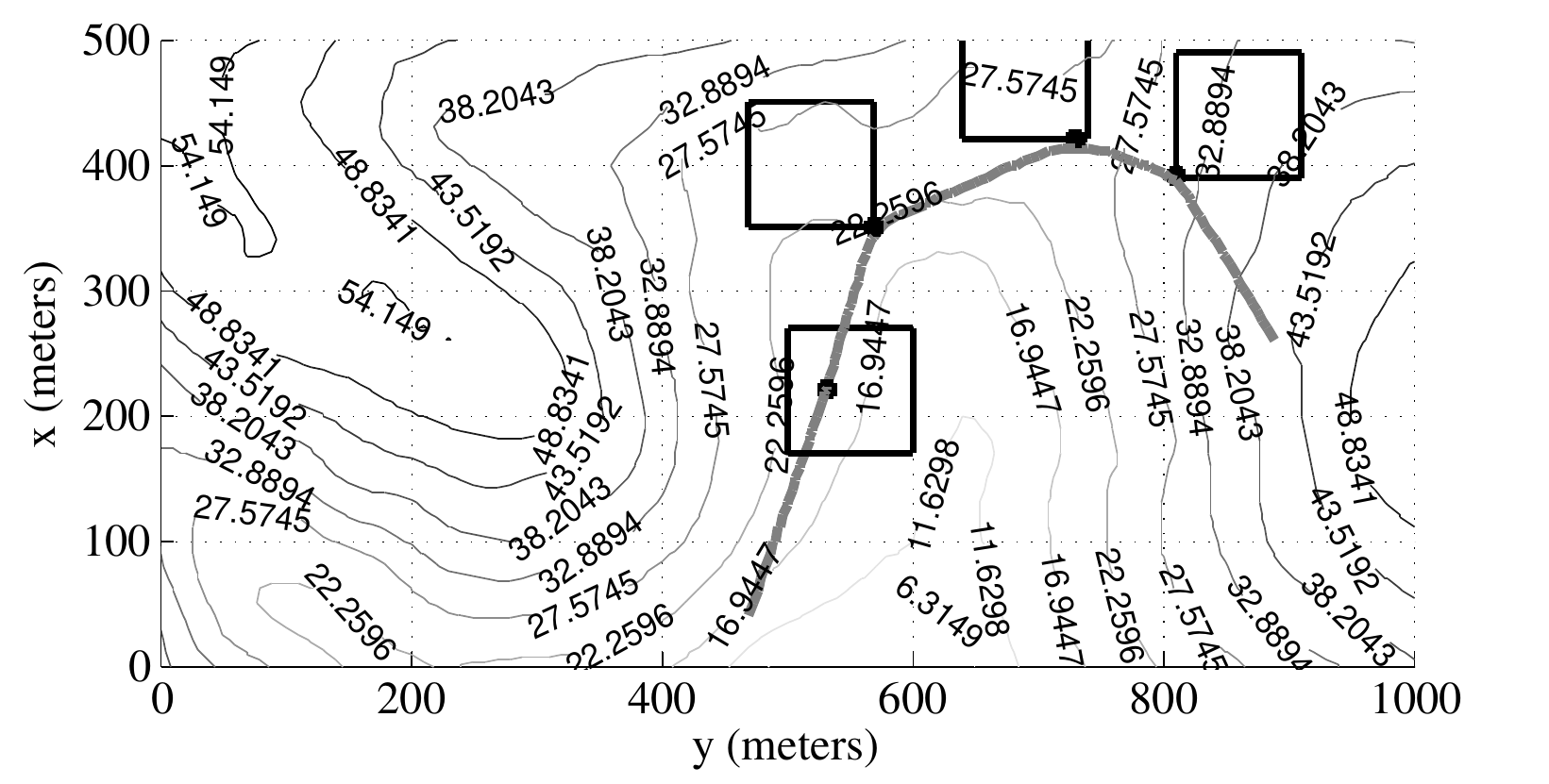}
	\caption{The horizontal alignment corresponding to the starred Pareto optimal point in Figure \ref{fig38}. The black squares represent the box constraints on the intersection points and the block points represent the intersection points.}
\end{figure}

\end{document}